\documentclass[11pt,reqno]{amsart}
\usepackage[margin=3.6cm]{geometry}

\usepackage{graphicx}
\usepackage{amssymb}
\usepackage{wrapfig}
\usepackage{bbm}

\newcommand{\A}{{\mathbb A}}
\newcommand{\R}{{\mathbb R}}
\newcommand{\N}{{\mathbb N}}

\newcommand{\C}{{\mathbb C}}
\newcommand{\Z}{{\mathbb Z}}

\newcommand{\ep}{\varepsilon}
 \newcommand{\repi} {\Re}
 \newcommand{\impi} {\Im}
  \newcommand{\Exp} {\mathrm{Exp}}

 \renewcommand{\Im} {\mathrm{Im\,}}
 \renewcommand{\Re} {\mathrm{Re\,}}
 \DeclareMathOperator{\inj}{inj_M}

\newtheorem{theorem}{Theorem}
\newtheorem{corollary}[theorem]{Corollary}
\newtheorem{proposition}[theorem]{Proposition}

\theoremstyle{definition}

\newtheorem{remark}{Remark}

\begin{document}

\title[Nodal sets in forbidden regions ]
{Nodal sets of Schr\"odinger eigenfunctions  in forbidden regions }

\author[Y. Canzani]{Yaiza Canzani}
\address[Y. Canzani]{ Institute for Advanced Study and  Harvard University.\medskip}
 \email{canzani@math.ias.edu}

\author[J. Toth]{John  A. Toth}
\address[J. Toth] {Department of Mathematics and Statistics, McGill University.}
 \email{jtoth@math.mcgill.ca}

\thanks{Y.C. was partially supported by an NSERC Postdoctoral Fellowship and by NSF grant DMS-1128155. J.T. was partially supported by NSERC grant  OGP0170280.}

\begin{abstract}
This note concerns the nodal sets of eigenfunctions of semiclassical Schr\"odinger operators acting on compact, smooth, Riemannian manifolds, with no boundary. We prove that if $H$ is a separating hypersurface that lies inside the classically forbidden region, then $H$ cannot persist as a component of the zero set of infinitely many eigenfunctions. In addition, on real analytic surfaces,  we obtain sharp upper bounds for the number of intersections of the zero sets of the Schr\"odinger eigenfunctions with a fixed curve  that lies inside the classically forbidden region.

\end{abstract}

\maketitle


Let $(M,g)$ be a smooth, compact, Riemannian manifold with no boundary. Write $\Delta_g$ for the Laplace operator, and given any smooth potential $V \in C^\infty(M;\R),$ consider the Schr\"odinger operator acting on $L^2(M)$ defined as 
$$P(h)=-h^2\Delta_g + V,$$
where $h \in (0,1]$.
Let $E \in \R$ be a regular value for the total energy function $p(x,\xi)=|\xi|^2_{g_x} + V(x)$ defined on $T^*M$, and  write $\Omega_E$ for the the classically forbidden region
$$\Omega_E:= \{ x \in M:\;   V(x)>E \}.$$
In this paper we study the nodal sets  of Schr\"odinger eigenfunctions (with energy close to $E$)  inside the classically forbidden region, in the semiclassical limit $h \to 0^+$.
 Consider $L^2$-normalized Schr\"odinger eigenfunctions  $\{\phi_h\}$ with
\begin{equation} \label{E: eigenfns}
P(h) \phi_h = E(h) \phi_h \qquad\text{and}\qquad E(h)=E+o(1) \;\; \text{as}\; \;h\to 0^+. 
\end{equation}

There is a large literature  devoted to the study of the zero sets of Laplace eigenfunctions,
 $$Z_{\phi_h} = \{ x \in M:\;  \phi_h(x)=0 \},$$
 on compact manifolds. We refer the reader to \cite{Z2} for a detailed list of references.  
The Hausdorff measure of the zero sets,  their distribution properties, the number of nodal domains and their inner radius, have been extensively studied (although many open problems remain, even for surfaces). More generally, it is natural to study the  properties  of  zero sets of Schr\"odinger eigenfunctions inside the classically allowed region where $V< E.$ Many of the known results in the homogeneous case where $V=0$ extend to  Schr\"{o}dinger eigenfunctions in the allowable region (see \cite{Jin}).
In contrast, very little is known about the zero sets of Schr\"odinger eigenfunctions inside the classically forbidden region where $V>E.$ 
 In dimension one, it is known that the eigenfunctions of the Harmonic Oscillator have no zeros in the forbidden region and in recent work, Hanin-Zelditch-Zhou  \cite{HZZ}  have proved that in any higher dimension the expected value of the measure of the  zero set of random eigenfunctions of the harmonic oscillator inside any ball is of order $h^{-1/2}$. 
We are not aware of any other results addressing the behavior of zero sets of Schr\"odinger eigenfunctions inside the classically forbidden region.

Our first result addresses the issue of {\em nodal persistence}: Can a fixed  hypersurface $H$ be contained in the nodal set of an infinite subsequence of eigenfunctions?  This question was answered on the flat torus  $\mathbb T^n$ by J.~Bourgain and Z.~Rudnick in \cite{BR}. They proved that if  $V=0$ and  $H\subset \mathbb T^n$ is a hypersurface with non-zero principal curvatures, then $H$ cannot lie within the zero set of infinitely many eigenfunctions.  
On general manifolds with  $V \neq 0,$ we prove that no embedded separating hypersurface contained entirely within the forbidden region $\Omega_E$ can persist as part of the zero set for infinitely many eigenfunctions.

\begin{figure}[h!]
\includegraphics[width=14cm]{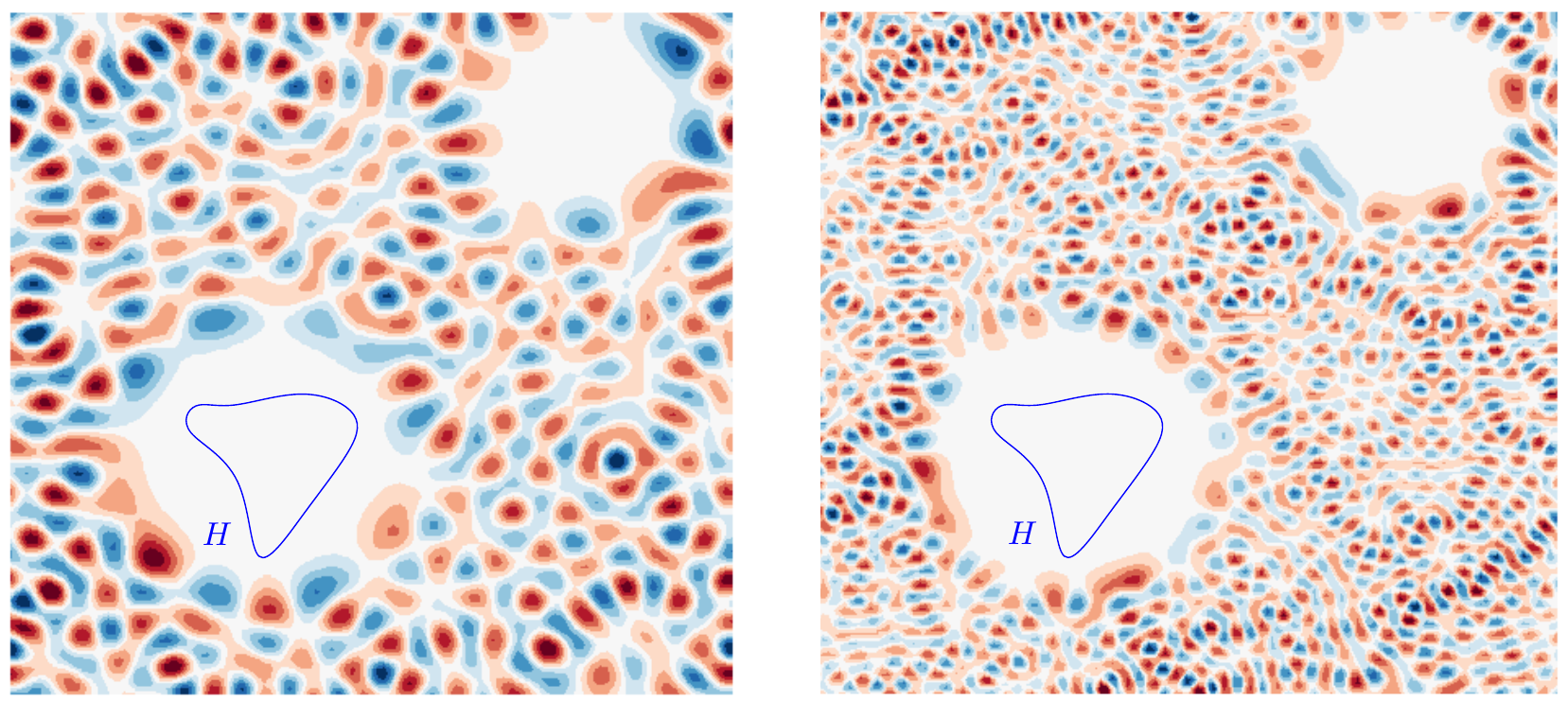}
\caption{\small Level sets of eigenfunctions of $-h^2\Delta_g+V$ on a square torus, where $E=1$ and $V$ is a periodized sum of two bumps $4 e^{-10((x+0.3)^2+(y+0.3)^2)} + 3 e^{-15((x-0.6)^2+(y-0.7)^2)}.$ The pictures correspond to $h=0.01$ and $h=0.005$ respectively. Tones of blue describe negative values, tones of red describe positive values. In the plot the value zero (and very small values as well) are depicted as white.}
\end{figure}
\vspace{-0.2cm}
\begin{theorem} \label{T: main 1}
Let $(M,g)$ be a smooth, compact,  Riemannian manifold with no boundary and let $V \in C^\infty(M)$. Consider a sequence  $\{\phi_{h}\}$ as in \eqref{E: eigenfns}.
Suppose that $H \subset \Omega_E$ is an  embedded separating hypersurface that encloses a bounded domain contained in $\Omega_E$.
 Then, there exist constants $C_H>0$ and $h_0 >0$ such that
 $$ \| \phi_{h} \|_{L^2(H)} \geq e^{- {C_H}/{h} } \qquad \text{and} \qquad  \|\partial_\nu \phi_{h} \|_{L^2(H)} \geq e^{- {C_H}/{h} },$$
 for all $h \in (0,h_0].$ 
In particular, for every  subsequence $\{h_{j}\}_j$ with $h_{j} \to 0$ as $j \to \infty,$ there is an integer $j_0 \geq 1$ with the property that for all $j \geq j_0$
$$H  \nsubseteq  Z_{\phi_{h_{j}}}.$$ 
\end{theorem}

\begin{remark} 
 Theorem \ref{T: main 1}  extends to   the case where $M= \R^2,$ provided that $V \in C^{\omega}(\R^2;\R)$  extends holomorphically  to a complex wedge domain $M^{\C}= \{ z \in \C^2: \,|\Im z| \leq \frac{1}{C} \langle z \rangle \}$ and that it satisfies $V(x) \geq C' \langle x\rangle^k$ for some $k \in \Z^+$  as $|x| \to \infty.$ Here, $C$ and $C'$ are positive constants.
 \end{remark}

Assume from now on that $(M,g)$ is a  compact, real analytic surface and let $H \subset \Omega_E$ be a  real analytic closed curve that bounds a region inside $\Omega_E$. 
Unique continuation results like the one in Theorem \ref{T: main 1} have important implications for the study of asymptotic oscillation properties of eigenfunctions, including estimates for the intersection number $ \#  \{ Z_{\phi_h} \cap H \}$ (see for example \cite{TZ,HT,CT}). Our second result is a deterministic upper bound for the nodal intersection with a fixed real-analytic curve $H$ on a Riemannian surface contained in the classically forbidden region.
\begin{theorem} \label{T: main 2}
Let $(M,g)$ be a  compact, real analytic surface with no boundary. Let $\{\phi_h\}$ be  real valued eigenfunctions satisfying \eqref{E: eigenfns}, where we also assume that the potential $V$ is real-analytic. Suppose that $H \subset \Omega_E$ is a simple, closed, real analytic curve that bounds a region inside $\Omega_E$. Then, there exists $C_H>0$ and $h_0>0$ such that 
$$  \#  \{ Z_{\phi_h} \cap H \}  \leq \frac{C_H}{ h},$$
for all $h \in (0, h_0]$.
\end{theorem}

To prove  Theorem \ref{T: main 2} we use the restriction lower bound in Theorem \ref{T: main 1} together with a potential layer formula for the eigenfunctions inside the forbidden region. Bounding the number of zeros on the fixed curve  is then reduced to estimating the  complexification of a particular Green's operator in the forbidden region.  We control the complexification of the Green's operator  using off-diagonal decay estimates for the real kernel  (see Proposition \ref{P: Greens operator}) together with  $h$-analytic Cauchy estimates recently proved by L. Jin in \cite{Jin}.

  For individual eigenfunctions, one can see that the $C_H h^{-1}$ bound in Theorem \ref{T: main 2} is sharp on surfaces of revolution (see Section \ref{S: example}) and agrees with the upper bound in Yau's conjecture [Y1,Y2] for nodal volume in the homogeneous case.
Nevertheless, it is reasonable to expect that in many cases one should be able to improve on this bound.
As mentioned, for random  eigenfunctions of isotropic harmonic oscillators, the computations in \cite{HZZ} show that the expected value of the nodal lengths  in the classically forbidden region are of order $h^{-1/2}$. Consequently, at least for random waves, it is reasonable to expect generic intersection bounds of the form
$ \#  \{ Z_{\phi_h} \cap H \}= O_H( h^{-1/2}) $ in the case for which the forbidden region is unbounded.
We hope to return to this question elsewhere.


\subsection{Organization of the paper}
In Section \ref{S: restriction bounds} we  prove Theorem \ref{T: main 1} using an elementary argument with Green's formula and quantitative unique continuation for the eigenfunctions. In Section \ref{S: nodal sets}, we study nodal intersection bounds by reproducing the eigenfunctions in the forbidden region using a suitable Green's operator whose complexification we need to control. Assuming that we have suitable bounds on the complexification of the Green's operator we then prove Theorem \ref{T: main 2} using the restriction lower bound in Theorem \ref{T: main 1}. In Section \ref{S: example}, we show that the upper bound in Theorem \ref{T: main 2} is sharp. 
In Section \ref{S: Greens operator} we give a detailed analysis of the kernel of  the  Green's operator on a compact manifold. In particular, we show that the kernel can be locally complexified away from the real diagonal  $\{(x,x) \in M \times M \}$ and obtain exponential decay estimates in $h$ for the complexified kernel.


\subsection{Acknowledgements} The authors would like to thank Christopher Wong for sharing his {\small MATLAB} code to compute estimations of Schr\"odinger eigenfunctions on the square torus.


\section{$L^2$-lower restriction bounds} \label{S: restriction bounds}

We note that because of the quantum tunnelling effect, the wave functions are known to have positive mass inside the classically forbidden region. Indeed, by Carleman type estimates \cite[Theorem 7.7]{Zw}, for every open set $U \subset \Omega_E$ there exists a positive constant $C = C(U)>0$ for which 
\begin{equation}\label{E: unique continuation}
\| \phi_{h} \|_{L^2(U)} \geq e^{- {C}/{h} }, \qquad \text{as}\;\; h\to 0^+. 
\end{equation}
The result in  Theorem \ref{T: main 1} is  an analogue of the lower bound in  \eqref{E: unique continuation} for the eigenfunction restricted to a hypersurface $H \subset \Omega_E$ and is a crucial step in the proof of Theorem \ref{T: main 2}. The only condition that we impose on $H$ is that it must bound a domain that is contained entirely inside $\Omega_E$.  

We note that the exponential lower bound in Theorem \ref{T: main 1} is quite delicate since despite the fact that inside the forbidden region  the eigenfunctions have positive mass, they are exponentially small in $h.$ Indeed, consider the Agmon metric $g_E = (V-E)_{+} g$ and associated distance function  $d_E(x):= \text{dist}_{g_E} (x, M \backslash\Omega_E).$ By the standard Agmon estimates \cite[Proposition 3.3.4]{Hel}, it follows that for any $\epsilon >0,$
\begin{equation} \label{E: Agmon}
| \partial_{x}^{\alpha} \phi_h(x) | = O_{\epsilon,\alpha}\Big( e^{ \frac{- d_{E}(x) + \epsilon}{h} } \Big) 
\end{equation}
locally uniformly in $x \in \Omega_{E}$.  In particular, given a smooth hypersurface $H \subset \Omega_E,$ it follows from \eqref{E: Agmon} that  for $d_E(H):=\min\{ d_E(x): \, x\in M\}$,  one has
\begin{equation} \label{E: upper}
\| \phi_h \|_{L^2(H)} = O\Big(e^{-\frac{d_E(H) + \epsilon}{h}} \Big).
 \end{equation}
We may then view Theorem \ref{T: main 1} as a partial converse to (\ref{E: upper})  under the assumption that $H \subset \Omega_E$ is a separating hypersurface.

\subsection{Proof of Theorem \ref{T: main 1}}
Let $H \subset \Omega_E$ be a separating hypersurface that bounds a smooth domain $M_H \subset \Omega_E.$
Since $M_H\subset \Omega_E$ and $ E(h)=E+o(1)$ as $h\to 0$, it follows that if $h_E>0$  is sufficiently small, 
 then there exists $C_E>0$ so that 
 \begin{equation}\label{E: V-E}
 V(x) -E(h) \geq C_E \qquad \text{  for all }\;\; x \in M_H
 \end{equation}
 and all $h \in (0,h_E]$.

 By Green's Theorem,
\[\int_{M_H} |h \nabla_g \phi_h|^2 dv_g + \int_{M_H} (h^2 \Delta_g) \phi_h  \, \overline{\phi_h} dv_g = h^2\int_{\partial M_H} \partial_\nu  \phi_h  \, \overline{\phi_h }d\sigma_g,\]
where $\nu$ is the outward normal vector and $\sigma_g$ is the induced volume measure on $\partial M_H$. 
Thus, since $-h^2\Delta_g \phi_h + V \phi_h= E(h) \phi_h$, it follows that 

\begin{align*}
 \| h \nabla_g \phi_h \|^2_{L^2(M_H)}  + \langle (V- E(h)) \phi_h, \phi_h \rangle_{L^2(M_H)}=h^2 \langle \partial_{\nu} \phi_h,  \phi_h \rangle_{L^2(\partial M_H)} .
\end{align*}

\bigskip

Using the non-negativity of $  \| h \nabla_g \phi_h \|^2_{L^2(M_H)}$ and \eqref{E: V-E} we obtain that for all  $h \in (0,h_E]$
\begin{align} \label{E: Green1}
 C_E\| \phi_h \|^2_{L^2(M_H)} \leq h^2 \langle \partial_{\nu} \phi_h ,  \phi_h \rangle_{L^2(\partial M_H)} .
\end{align}
An application of the   Cauchy-Schwarz inequality in \eqref{E: Green1} gives
\begin{align} \label{E: Green 2}
  C_E h^{-2} \| \phi_h \|_{L^2(M_H)}^2 \leq \| \phi_h \|_{L^2(H)} \, \| \partial_{\nu}\phi_h \|_{L^2(H)}.
\end{align}
By the unique continuation lower bound \eqref{E: unique continuation} there exists $C>0$ with  $ \| \phi_h \|_{L^2(M_H)}^2 \geq e^{-C/h}$ for $h$ small enough and therefore there exists $h_E>0$ such  that 
   \begin{equation*}
 C_E h^{-2}  e^{-\frac{C}{h}}  \leq  \| \phi_h \|_{L^2(H)} \| \partial_{\nu}\phi_h \|_{L^2(H)},
\end{equation*}
 for all  $h \in (0,h_E]$.   Theorem \ref{T: main 1}  then follows from  the Agmon estimates in \eqref{E: Agmon}.\\
 \qed
 


\section{Nodal intersection bounds}\label{S: nodal sets}

Here we present the proof of Theorem \ref{T: main 2} (Section \ref{S: Proof of Thm 2}), and show that the upper bound on the number of nodal intersections is saturated for surfaces of revolution (Section \ref{S: example}).

\subsection{Proof of Theorem \ref{T: main 2}}\label{S: Proof of Thm 2}
We continue to assume that $H \subset \Omega_E$, but here we make the additional assumption that the Riemannian manifold $(M,g),$ the potential  $V$ and the hypersurface $H$ are all real-analytic.  Also, in the following, we are only interested in  the case in which $M$ is a surface. Let  $q:[0,2\pi] \to H$  be a $C^{\omega}$, $2\pi$-periodic, parametrization of $H$. 
To bound the number of zeros of $\phi_h \circ q: [0, 2\pi] \to \R$ we consider its holomorphic extension $(\phi_h \circ q)^\C:H_\tau^\C \to \C$ to the complex strip  
  $$H_{\tau}^{\C} = \{ t \in \C: \; \Re t \in [0,2\pi],  \; |\Im t| < \tau \}$$
  for some $\tau>0$,
  and use that $  \#\{Z_{\phi_h} \cap H \} \leq  \# \{ Z_{(\phi_h \circ q)^\C} \cap H_{\tau}^{\C}\}$. Then, the zeros of $(\phi_h \circ q)^\C$ are studied using the Poincar\'e-Lelong formula:
  $$\partial \overline{\partial} \log |(\phi_h \circ q)^\C(z)|^2=\sum_{z_k \in  Z_{(\phi_h \circ q)^\C} } \delta_{z_k}(z).$$
According to \cite[Proposition 10]{TZ},  there exists $C>0$ so that
 \begin{equation} \label{E: [TZ]}
  \#\{Z_{\phi_h} \cap H \} \leq \# \{ Z_{(\phi_h \circ q)^\C} \cap H_{\tau}^{\C}\} \leq  C \, \max_{t \in H_{\tau}^{\C}} \log |F_{h}^{\C}(t)|,
   \end{equation}
where $F_h^{\C}(t)$ with $t \in H_{\tau}^{\C}$ is the holomorphic continuation of  the normalized eigenfunction traces 
 \begin{equation}\label{E: U(t)}
  F_{h}(t):= \frac{\phi_h(q(t))}{\|\phi_h \|_{L^2(H)}}.
  \end{equation}
Note that by Theorem \ref{T: main 1} we know that $\|\phi_h \|_{H} > e^{-C_H/h}$ for $h\in (0, h_0]$ with $h_0$ sufficiently small, and this implies that $F_h(t)$ is well defined. 
 
 It follows that we shall need to control the complexification $F_h^{\C}(t)$ to obtain upper bounds on $\#\{Z_{\phi_h} \cap H \} $. 
 Without loss of generality we assume that $H \subset \, \text{int}  (\Omega_H)$ where  $\Omega_H \subsetneq \Omega_E$ is  a domain  containing $\gamma$,  a closed $C^\omega$ curve,  as  its boundary. 
Our goal is to find a double layer jumps formula that reproduces $\phi_h(x)$ for $x \in H$ in terms of its values along $\gamma$.

Let $\chi \in C^{\infty}_{0}(M, [0,1])$ with $\chi (x) = 1$ for all $x \in \Omega_H$ and with $\text{supp} \chi \subset \Omega_E$. Consider the auxiliary global metric given by
\begin{equation} \label{jacobi}
g_{_{\Omega_E,h}}(x) := (V(x) - E(h))  \chi(x) g(x) + (1-\chi(x)) g(x), \qquad x \in M.
 \end{equation}
 From now on, to simplify notation, we simply write $g_{_{\Omega_E}}$ for $g_{_{\Omega_E,h}}$ since the dependence of the latter on $h$ is only in the constant eigenvalue term and is of no real  consequence as far as $h$-pseudodifferential calculus is concerned.

Then, since $n=2,$ it follows that $\Delta_{g_{_{\Omega_E}}} =  [ (V-E(h)) \chi + (1-\chi) ]^{-1} \, \Delta_g$ and so, in particular,  for all $x \in \Omega_H$,
\begin{equation} \label{eigenequation}
( -h^2 \Delta_{g_{_{\Omega_E}}}  + 1) \phi_h(x) = -h^2 (V(x) - E(h))^{-1} \Delta_g \phi_h(x) + \phi_h(x) = 0.  \end{equation}

 We consider the Green's operator
 $$G(h) = ( -h^2 \Delta_{g_{_{\Omega_E}}} + 1)^{-1}.$$ 
Since $(-h^2 \Delta_{g_{_{\Omega_E}},y} + 1) G(x,y,h)=\delta_x(y)$,  (\ref{eigenequation}) implies that for all $x \in M_H$
\begin{align} \label{green}
&\phi_h(x) = \notag \\
&= \int_{\Omega_H} (-h^2 \Delta_{g_{_{\Omega_E}},y} + 1) G(x,y,h)  \phi_h(y) dv(y)
 - \int_{\Omega_H} G(x,y,h)   (-h^2 \Delta_{g_{_{\Omega_E}} } + 1) \phi_h(y)  dv(y). 
\end{align}
By Green's formula, it then follows that for $x \in \Omega_H,$
\begin{align} \label{potential layer}
\phi_h(x) =  h^2 \Big(  \int_{\gamma} G(x,y,h)  \partial_{\nu_y} \phi_h(y) d\sigma(y)-\int_{\gamma} \partial_{\nu_y} G(x,y,h) \phi_h(y) d\sigma(y)  \Big), \end{align}
where $\nu_y$ is the outward normal vector at $y \in \gamma$.
Let   $r:[0,2\pi] \to \gamma$ be a $C^{\omega}$ parametrization of $\gamma$. Restriction of the outgoing variable $x$ in (\ref{potential layer}) to $H$ yields the potential layer  formula
\begin{equation} \label{potential layer 2}
\phi_h(q(t)) = h^2 \int_{\gamma} G(q(t),r(s),h)   \, \partial_{\nu_y} \phi_h(r(s)) d\sigma(s)- h^2\int_{\gamma} \partial_{\nu_y} G(q(t),r(s),h) \phi_h(r(s)) d\sigma(s) .
\end{equation}

In order to control $F_h^\C(t)$ in \eqref{E: [TZ]},  one needs an upper bound for the holomorphic continuation $(\phi_h \circ q)^\C.$  The latter amounts to  estimating the complexification of  \eqref{potential layer 2}.

  Given $\tau>0$, let $q^{\C}(t)$ denote the holomorphic continuation of the $C^{\omega}$ parametrization $q:[0,2\pi] \to H$ to the strip $H_\tau^\C$. We claim the following result.
 
  \begin{proposition}\label{P: complexification}
 Suppose $d_{g}(H, H_\delta)>\ep$ for some $\ep >0.$  Then, there exist  constants $C(\ep)>0$, $\tau(\ep)>0$, $h_0(\ep)$, such that for $h \in (0,h_0(\ep)],$
 \begin{equation} 
 |G^{\C}(q^{\C}(t),r(s),h)| = O(e^{-\frac{C(\ep)}{h}}) \qquad \text{and} \qquad |\partial_{\nu_y}G^{\C}(q^{\C}(t),r(s),h)| = O(e^{-\frac{C(\ep)}{h}}),  \end{equation}
 uniformly for $(t,s) \in H_{\tau(\ep)}^{\C} \times [0,2\pi].$  
 \end{proposition}
 
Proposition \ref{P: complexification} is a consequence of a more general result, Theorem \ref{T: Greens operator complexification}, which we prove   in Section \ref{S: Greens operator} (see Remark \ref{link}). In Theorem \ref{T: Greens operator complexification}, we show that the kernel of the Green's operator $G(h)$ along with its derivatives can  be locally complexified off-diagonal maintaining the exponential decay exhibited in the real domain.

  Substitution of the estimates in Proposition \ref{P: complexification}  into the complexification of \eqref{potential layer 2}, combined with an application of  Theorem \ref{T: main 1}   and the Cauchy-Schwarz inequality gives the existence of positive constants $C, h_0$ and $d_H$ such that
 \begin{equation} \label{forbidden}
 |F_h^{\C}(t)| \leq C e^{-C_1(\ep)/h} \left( \frac{ \| \phi_h \|_{L^2(\gamma)} }{ \|\phi_h \|_{L^2(H)} } + \frac{ \| \partial_{\nu} \phi_h \|_{L^2(\gamma)} }{ \|\phi_h \|_{L^2(H)} } \right) = O(e^{d_H/h}) 
 \end{equation}
 for all $h \in (0, h_0]$. Then, by \eqref{E: [TZ]}, there exists $C_H>0$ such that 
 \begin{equation}
  \#\{Z_{\phi_h} \cap H \} \leq C_H h^{-1},
  \end{equation}
as desired. Assuming that Proposition \ref{P: complexification} holds, this concludes the proof of  Theorem \ref{T: main 2}. \qed

\subsection{Estimates in classically allowable versus forbidden regions}  It is interesting to contrast the growth estimates in (\ref{forbidden}) with the case where $H$ is contained in the classically allowable region. For example, when $H \subset \Omega$ where $\Omega \subset \R^2$ is a piecewise-analytic planar domain, and $\phi_h$ is a homogeneous eigenfunction (satisfying either Neumann or Dirichet boundary conditions), one can show that (see \cite{TZ} Lemma 11),
\begin{equation} \label{allowable}
 |F_h^{\C}(t)| \leq C e^{C_2/h} \left( \frac{ \| \phi_h \|_{L^2(\partial \Omega)} }{ \|\phi_h \|_{L^2(H)} } + \frac{ \| \partial_{\nu} \phi_h \|_{L^2(\partial \Omega)} }{ \|\phi_h \|_{L^2(H)} } \right) = O(e^{ d_H'/h}). 
 \end{equation}
 In (\ref{allowable}), the constant $C_2 = \max_{(q,r) \in H \times \partial \Omega} \Re i d^{\C}(q^{\C}(t),r(s))$ where $d^{\C}$ is the complexified distance between $\partial \Omega $ and $H$, and $q,r$ are as defined in Section \ref{S: Proof of Thm 2}. Thus, $C_2$ is {\em positive} in contrast with the {\em negative} constant $-C_1(\ep)$ appearing in the forbidden case. This is due to the fact that the Green's kernel $G(q(t),r(s),h)$ is a semiclassical pseudodifferential operator that decays exponentially off the diagonal and so does the corresponding local complexification $G^{\C}(q^{\C}(t),r(s),h)$ (see Theorem \ref{T: Greens operator complexification}). In the allowable region, the Green's kernel $G(q(t),r(s),h)$ is replaced with the restriction of the free Helmholtz Green's kernel $G_{\R^2}(q(t),r(s),h)$ in $\R^{2}$ which has the WKB asymptotics
 $G_{\R^2}(q(t),r(s),h) \sim_{h \to 0^+} (2\pi h)^{1/2} e^{id(q(t),r(s))/h} (a_0(t,s) + a_1(t,s) h + \cdots),$ provided $\inf_{(q,r) \in H \times \partial \Omega} d(q,r) >0.$
 This is the kernel of an $h$-Fourier integral operator  and the phase factor $e^{id/h}$ blows up exponentially in $h \to 0^{+}$ upon complexification in $q(t)$, unlike in the $h$-pseudodifferential case where there is off-diagonal exponential decay in $h.$ In view of the Jensen-type growth estimate in (\ref{E: [TZ]}), it follows that the {\em constant} in the $O_H(h^{-1})$ intersection bound in the forbidden region is smaller than the one for the allowable region. However, as the next example shows, the $h^{-1}$-rate cannot be improved in general. We hope to return to discuss these issues in more detail elsewhere.


\subsection{The example of a convex surface of revolution} \label{S: example}
Here we show that the upper bound in Theorem \ref{T: main 2} is sharp.
 To do this,  consider a convex surface of revolution generated by rotating a curve $y=f(r)$ about the $r$-axis with  $f\in C^w([-1,1], \R),$  $f(1)=f(-1)=0$ and in addition require that    $f''(r)<0$ for all $r \in [-1,1]$ so that the surface is strictly convex. Let $M$ be the corresponding surface of revolution parametrized by
$$[-1,1] \times  [0, 2\pi)  \to \R^3,$$
$$(r, \theta) \mapsto (r\,,\, f(r) \cos(\theta)\,,\, f(r) \sin(\theta)).$$
In these coordinates $M$ inherits a Riemannian metric $g$ given by 
$$g=  w^2(r) dr^2 + f^2(r) d\theta^2$$
where we have set $w(r):=\sqrt{1+ (f'(r)) ^2}$. 
Consider on $(M,g)$ an analytic potential $V(r,\theta)=V(r)$  independent of the angular variable, and let $E$ be a regular energy level for the hamiltonian corresponding to $-h^2\Delta_g +V$ with $ \min V < E < \max V.$ We shall construct a curve $H$  contained in the forbidden region $\{V>E\}$  and a sequence of real valued solutions $\{\phi_{h_k}\}_k$ of 
$(-h^2 \Delta_g + V )\phi_{h_k} = E(h) \phi_{h_k}$
where $E(h_k)=E +o(1)$ as $k \to +\infty$, so that
$$\#\{\phi_{h_k}^{-1}(0) \cap H\} \geq 2h_k^{-1}.$$
Consequently, the  $O(h^{-1})$ bound in Theorem \ref{T: main 2} is sharp.
  \begin{figure}[h!]
\includegraphics[height=3.5cm]{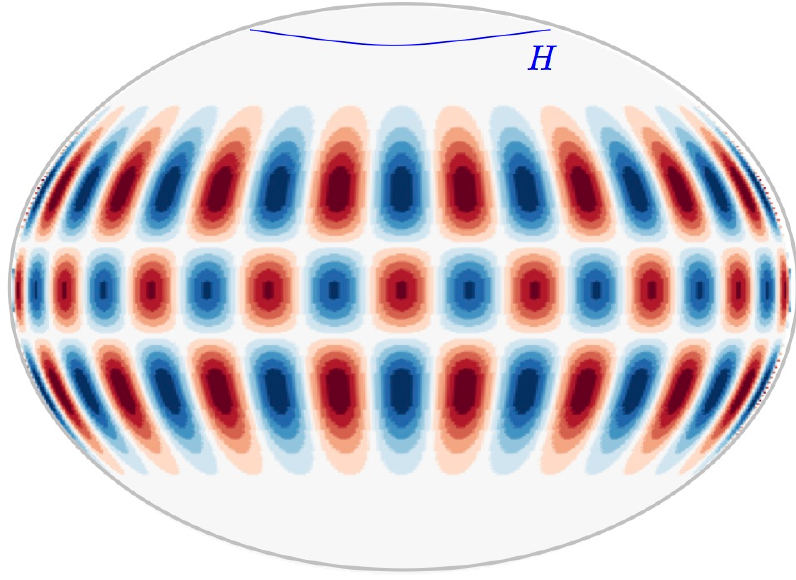}
\caption{Illustration of the level sets of an eigenfunction in an ellipsoid where the potential is concentrated close to the poles. The colouring scheme is the same as in Figure 1.}
 \end{figure}

We seek eigenfunctions of the form $\phi_h(r, \theta)= v_h(r) \psi_h(\theta)$ that solve
$$(-h^2 \Delta_g + V )\phi_h = E(h) \phi_h.$$
Since the Laplace operator in the coordinates $(r,\theta)$ takes the form
\begin{equation*}
\Delta_g=\frac{1}{w(r)f(r)} \frac{\partial} {\partial r} \left( \frac{f(r)}{w(r)} \frac{\partial} {\partial r} \right)  + \frac{w^2(r)}{f^2(r)}  \frac{\partial^2} {\partial \theta^2},
\end{equation*}
we have that  the functions $v_h$ and $\psi_h$ must satisfy
\begin{equation}\label{psi}
-  \frac{d^2}{d\theta^2}\psi_h(\theta) = m_h^2\,  \psi_h(\theta)
\end{equation} 
and
\begin{equation}\label{v}
-h^2 \frac{f(r)}{w(r)} \frac{d} {d r} \left( \frac{f(r)}{w(r)} \frac{d} {d r} v_h(r) \right) + f^2(r)(V(r)-E(h))v_h(r)=  - m_h^2 \,  h^2 w^2(r) v_h(r), 
\end{equation}
for some $m_h \in \mathbb Z$. From now on let $\{h_k\}_k$ be a decreasing sequence with $h_k \to 0^+$ as $k \to +\infty$ and    such that  $m_{h_k}= 1/h_k \in \Z$. One can choose the solution of (\ref{psi}) to be 
\begin{equation}\label{psi solution}
\psi_{h_k}(\theta)=e^{{i \theta/h_k} }.
\end{equation}
To deal with the radial part $v_{h_k},$ one makes the change of variables $s \mapsto r(s) = \int_{0}^{s} \frac{f(\tau)}{w(\tau)} d\tau$ in \eqref{v} and it follows that 
\begin{equation}\label{v 2}
\left(-h_k^2 \frac{d^2} {d s^2}+ f^2(r(s))(V(r(s))-E(h))+ w^2(r(s)) \right) v_{h_k}(r(s))=0.
\end{equation}

To finish the argument,  let $r_0 \in \{r \in [-1,1]:\; V(r) >E\}$ and set $H_{r_0}$ to be the curve
$$H_{r_0}:=\{(r, \theta):\; r=r_0\}.$$
Since $\cos(\theta/h_k)$ has $2h_k^{-1}$ zeros for $\theta \in [0, 2\pi)$, our claim is established once we set
$$\phi_{h_k}(r, \theta) := \Re(e^{i\theta/h_k}\;v_{h_k}(r))$$
with $v_{h_k}$ a solution to \eqref{v 2}.

\section{The Green's operator $G(h)$ and its complexification}  \label{S: Greens operator}
Let $(M,g)$ be a compact, real analytic Riemannian manifold of dimension $n.$ We   consider here  the  associated Green's operator  $$G(h)= (-h^2 \Delta_{g} +1)^{-1}: C^{\infty}(M) \to C^{\infty}(M).$$ The purpose of this section is to study the complexification of the kernel $G(x,y,h)$ in the outgoing variable $x$. Before we state our main result  (Theorem \ref{T: Greens operator complexification}), we briefly review some of the complex analytic geometry that is needed in the formulation and proof of Theorem \ref{T: Greens operator complexification}.

\subsection{Grauert tube complexification of $M$}
By a theorem of Bruhat-Whitney,  $M$ has a unique complexification $M^\C$ with $M \subset M^{\C}$ totally real that generalizes the complexification of $\R^n$ to $\C^n$. 
The open Grauert tube of radius $\varepsilon$ is defined to be  
$$M_\tau^\C=\{ z \in M^\C: \; \sqrt{\rho_g}(z) < \tau\},$$
where $\sqrt {\rho_g}$ on $M^\C$ is the unique solution to the complex Monge-Ampere equation.
 For example, in the simplest model case when $M = \R^2$ and $M^{\C} = \C^2,$ one has $\sqrt \rho_g(z) =2 | \Im z|.$ 
 There is a maximal $\tau_{\max}>0$ for which $M_\tau^\C$ is defined \cite[Thm 1.5]{LGS}, and  $M_\tau^\C$ is a strictly pseudoconvex domain in $M_{\tau_{\max}}
^\C$ for all $\tau \leq \tau_{\max}$. 

For all $\tau \leq \tau_{\max}$, we identify the radius $\tau$ ball bundle $(BM)_{\tau} \subset TM$ with $(B^*M)_{\tau} \subset T^*M$ using the Riemannian metric. For $x \in M$ and $0<r<\text{inj}(M,g)$, we let $\exp_x:B_x(0,r) \to M$ be the geodesic exponential map defined on the geodesic ball $B_x(0,r) \subset T^*_x M $. We denote the lifted exponential map to all of $(B^*M)_{\tau} $ by
$$\Exp:(B^*M)_{\tau}  \to M, \,\,\,\, \Exp(x,\xi) =  \exp_{x}(\xi).$$ Since $(M,g)$ is real-analytic, for fixed $x \in M$ and $0<r<\inj$,  the geodesic exponential map $\exp_x: B_x(0,r)  \to M$   admits a holomorphic continuation $\exp_x^{\C}:(B_x(0,r))^{\C} \to M^\C$ in the fiber $\xi$-variables.
For $0<\tau< \tau_{\max}$, we define the associated complexified lifted map by
$$\Exp^\C: (B^*M)_{\tau} \to M^\C, \,\,\,\, \Exp^\C(x,\xi)= \exp_x^{\C} (i \xi).$$
The map $\Exp^\C$ gives a diffeomorphism between $(B^*M)_{\tau}  $ and $M_\tau^\C$ with the property that
$({\Exp^\C})^* (\rho_g)=|\cdot|_g.$  Consequently, $(B^*M)_{\tau} \cong M^\C_{\tau}$ as complex manifolds via $\Exp^{\C}$. Also, the map
\begin{equation}\label{pi_M}
 \pi_{_M}:\, M_\tau^\C \to M, \,\qquad \,\, \pi_{_M}(\Exp^\C(x, \xi))=x
 \end{equation}
is an analytic fibration. The fibers $\pi_{_M}^{-1}(M)$ correspond to imaginary directions over the totally real submanifold $M \subset M_\tau^\C.$

In general, we shall denote by $\phi^\C:M^\C_\tau \to \C$ the complexification of a real analytic function $\phi: M \to \C$.

Fix $x_0 \in M$. The map 
\begin{align*}
B_{x_0}(0,r) &\to M\\
\eta=r(x)  &\mapsto \exp_{x_0}(\eta) =x,
\end{align*}
is real analytic near the origin and so it can be holomorphically extended as
\begin{align*}
(B_{x_0}(0,r))^\C &\to M_\tau^\C \approx (B^*M)_{\tau}\\
 \eta + i \zeta= f(x,\xi) & \mapsto \exp_{x_0}^\C(\eta + i\zeta)= (x,\xi).
 \end{align*}
By Lemma 1.18 in \cite{LGS}, this coordinate system satisfies $f(x,0)=r(x)$ and $f(x_0,\xi)=i \xi$. 
Identifying the point $(x,\xi) \in B_\tau^*M$ with $\exp_x^\C(i\xi) \in M^\C_\tau$ as described above, one has $\pi_{_M}(x,\xi)=\pi_{_M}(\Exp^\C(x,\xi))=x=\exp_{x_0}(\eta)$. In view of Lemma 1.18 of \cite{LGS} we will use holomorphic coordinates $(\eta, \zeta)=f(x,\xi)$ on the complex manifold $(B^*M)_\tau$.

From now on, in a coordinate neighbourhood of $x_0 \in M,$ we write $z= z(x,\xi) \in M^\C_\tau$ for complex coordinates where \[z=(\Re z, \Im z)\]  with
\begin{equation}\label{normal coordinates}
 \repi z := \Re f(x,\xi), \quad \text{and}  \quad \impi z:= \Im f(x,\xi).
\end{equation}
Note  that with this notation   $\pi_{_M}(z)=\pi_{_M}(x,\xi)=\exp_{x_0}(\eta)=\exp_{x_0}(\Re f(x,\xi))=\exp_{x_0}(\Re z)$, and so $\pi_M (z)$ is identified with $\Re z$.

We also have that $\rho_g^{2} \in C^{\infty}(M^{\C}_{\tau})$ is a strictly plurisubharmonic exhaustion function, and by Taylor expansion around $\Im z =0$ it follows that $\rho_g^{2}(z) = 4 |\Im z|^2 + O(|\Im z|^{3}).$ Thus, for $z \in M_\tau^\C$ the function $\sqrt{\rho_g}(z)  \sim 2 |\Im z|$ and  it will sometimes be convenient to work on the subdomain of  $M_\tau^\C$ given by
$$ \Big\{z \in M_{\tau_{max}}^\C:\;  |\Im z|  \leq \frac{\tau}{C} \Big\}$$
with $C>0$ sufficiently large.

\subsection{Statement of the main result}
 Let $\ep>0$ and consider  the $\ep$-diagonal neighborhood $\Delta(\epsilon) = \{ (x,y) \in M\times M:  d_g(x,y) < \epsilon \}.$  Under the analyticity assumption on $(M,g),$ 
\begin{equation} \label{analytichypoelliptic}
G(\cdot, \cdot,h) \in C^{\omega}(M \times M \backslash \Delta(\epsilon)). \end{equation}
 Indeed, for   any $h \in (0,1),$   $-h^2\Delta_{g} + 1$ is a real-analytic partial differential operator that is uniformly elliptic in $h \in (0,1)$ and   for  fixed $y \in M$ and any $\epsilon >0,$
 \begin{equation} \label{analytic equation}
 ( -h^2\Delta_{g,x} + 1 ) G(x,y,h) = h^{-2} \delta(x-y) = 0 
  \end{equation}
 whenever $(x,y)\in (M \times M)\backslash  \Delta(\epsilon)$.   Consequently, (\ref{analytichypoelliptic}) follows by analytic hypoellipticity.
We may then consider the complexification $G^{\C}(z,y,h)$  in the outgoing variable  for $(z,y) \in M_{\tau(\ep)}^{\C} \times M$ with $d_g(\pi_{_M} z, y)>\ep$, where $\tau(\ep) >0$ is sufficiently small depending only on $\ep$.
The purpose of this section is to prove the following  asymptotic (in $h$) supremum bound for the holomorphic continuation $G^{\C}(z,y,h)$ in a small complex Grauert tube over the real off-diagonal domain.

 \begin{theorem}\label{T: Greens operator complexification}
Given  $\ep>0$  there exists a constant $\tau(\ep)>0$   such  that the Greens kernel admits a holomorphic extension $ G^\C(z,y,h)$ to $(z,y) \in M_{\tau(\ep)}^{\C} \times M$ with $d_g(\pi_{_M} z, y)>\ep$. Moreover, for any fixed $N_0 \in \N^n$, and $\alpha \in \N^n$ with $|\alpha| \leq N_0,$ there exists $ C=C(\ep,N_0)>0$ such that as $h \to 0^+$
   \begin{equation} \label{E: G complex Theorem}
 |\partial_y^\alpha G^{\C}(z,y,h)| = O\left(e^{-\frac{C}{h}}\right), 
  \end{equation}
 uniformly for $(z,y) \in M_{\tau(\ep)}^{\C} \times M$ with $d_g(\pi_{_M}z, y)>\ep$.
 \end{theorem}
 
 We prove  Theorem \ref{T: Greens operator complexification} in Section \ref{S: Greens kernel complexification}.
 
 \begin{remark}[Proof of  Proposition \ref{P: complexification}] \label{link} 
  Proposition \ref{P: complexification} follows directly from Theorem \ref{T: Greens operator complexification} as a special case. Indeed,
 since $d_g(H, \gamma) > \ep$ for some $\ep>0$,  by choosing the Grauert tube radius $\tau(\ep) >0$ small enough, the estimate in (\ref{E: G complex Theorem})  is satisfied by the Green's kernel associated with the extended Agmon metric $g_{_{\Omega_E}}$ defined in \eqref{jacobi}. We note that although $g_{_{\Omega_E}}$ is only globally $C^{\infty}$ on $M$, it is real-analytic in the forbidden region $\Omega_E.$ Thus, the Schwartz kernel, $G(x,y,h),$ of $G(h) = (-h^2 \Delta_{g_{\Omega_E} } +1)^{-1}$ is real-analytic for $(x,y) \in \Omega_E \times \Omega_E.$
  \end{remark}

To prove Theorem \ref{T: Greens operator complexification}, we must first  describe the  real kernel $G(x,y,h)$ in detail. In particular, we  prove exponential  off-diagonal decay estimates for the real kernel $G(x,y,h)$ in  Proposition \ref{P: Greens operator} below. Although such results are known to experts, we could not find a reference in the literature containing all the details we need here. As a result, for completeness, we carry out in detail the $h$-analytic parametrix construction in \cite{Sj} using the method of analytic stationary phase, keeping careful track of the various remainder terms. 

Consider the cut-off function
\begin{equation}\label{E: chi}
{\chi}(x,y) := \chi (d_g(x,y)),
\end{equation}
for $ \chi \in C^{\infty}_{0}([0,+\infty] ;[0,1])$, with   $\text{supp}\,\chi \subset [0,\text{inj} (M,g)]$ and
$$ \chi(t)=1 \quad \text{for} \quad  t\in \left[0, {\text{inj}(M,g)}/{2}\right].$$
We cover the manifold with coordinate patches so that if $\chi(x,y)\neq 0$ then $x$ and $y$ belong to a common coordinate neighborhood.

 \begin{proposition}\label{P: Greens operator}
 The Schwartz kernel of the Greens operator admits a decomposition of the form 
  \begin{align}\label{E: G theorem}
 G(x,y,h)= A_G(x,y,h) + R_G(x,y,h).
 \end{align}
Here,
 $$A_G(x,y,h)=\frac{\chi(x,y)}{(2\pi h)^n} \int_{\R^n}e^{ \frac{i}{h}   \langle g^{-1}_y \exp_{y}^{-1} (x),   \eta \rangle -\frac{1}{4h} d_g^2(x,y) \langle   g^{-1/2}_y \eta \rangle} \,a_G( x,y,\eta, h) \, d\eta,$$
with
$a_G( x,y,\eta, h)=\sum_{k=0}^{\frac{1}{C_0}[\frac{1}{h}]-1} h^k w_k(x,y,\eta) \in S^{0,-2}_{cl}$ for some $C_0>0$.
Also,  for each $\alpha, \beta \in \N^n$ there exists $C_{\alpha,\beta}>0$  and $h_0 = h_0(\alpha,\beta)>0$ such that  for all $h \in (0,h_0],$
\begin{equation}\label{E: R_G theorem}
 |\partial^\alpha_x \partial^\beta_y R_G(x,y,h)| = O(e^{-C_{\alpha,\beta}/h}),
 \end{equation}
uniformly in $x,y \in M.$
\end{proposition}

\begin{remark} \label{max}
We note that writing $G_y(x;h):= G(x,y,h)$ it follows that for $x \in M \backslash B_y(\ep)$ the function $G_y$ satisfies 
$P(h) G_y(x) = (-h^2 \Delta_g +1) G_y(x) = 0$, and since $P(h)$ is uniformly elliptic for $h < h_0,$ it follows from the maximum principle that
$$ \max_{ \{x:\, d_g(x,y) > \frac{1}{2} \text{inj}(M,g) \} } |G_y(x)| \leq  \max_{ \{x:\, d_g(x,y) =\frac{1}{2} \text{inj}(M,g) \} } |G_y(x)| = O(e^{-C/h}),$$
where the last estimate follows directly from Proposition \ref{P: Greens operator}. Thus, it suffices to bound $\chi(x,y) G(x,y,h)$ since the far off-diagonal part of the Green's kernel $(1-\chi(x,y)) G(x,y,h)$ is controlled by the former and is absorbed into the remainder term $R_G(x,y,h)$ in Proposition \ref{P: Greens operator}.
\end{remark}

We prove Proposition \ref{P: Greens operator} in Section \ref{S: Greens kernel}.
 In terms of normal coordinates centered at $y,$ 
\[A_G(x,y,h)= \frac{\chi(x,y)}{(2\pi h)^n}  \int_{\R^n}e^{ \frac{i}{h}   \langle x-y,   \eta \rangle -\frac{1}{4h} |x-y|^2 \langle\eta\rangle} \,  a_G(x,y,\eta,h) \,d\eta,\]

and so, restriction of $A_G(x,y,h)$ to points $x,y \in M$ with $d_g(x,y)>\ep >0$ gives
$$ \Re \Big( i   \langle g^{-1}_y \exp_{y}^{-1} (x),   \eta \rangle -\frac{1}{4} d_g^2(x,y) \langle   g^{-1/2}_y \eta \rangle \Big) \leq - \frac{C \ep^2}{4}  \langle \eta \rangle $$
for some $C>0$. As a consequence of Proposition \ref{P: Greens operator}, one gets the off-diagonal exponential decay estimates for $G(x,y,h).$

\begin{corollary}\label{C: exponential decay real kernel}
Let  $\ep>0$ and $N_0 \in \mathbb N$.  Then, there exists $C =C(\ep,N_0)>0$ such that for all $\alpha,\beta \in \N^n$ with $|\alpha| \leq N_0$ and $|\beta| \leq N_0,$ 
   \begin{equation} \label{E: G real Theorem}
 |\partial^\alpha_x \partial_y^\beta G(x,y,h)| = O \left(e^{-\frac{C}{h}}\right)
  \end{equation}
  as $h \to 0^+$, uniformly for $(x,y) \in M \times M$ with $d_g(x, y)>\ep$.
\end{corollary}

In order to avoid breaking the exposition at this point, we defer the proof of   Proposition \ref{P: Greens operator}    to Section \ref{S: Greens kernel} and proceed with the proof of Theorem \ref{T: Greens operator complexification}.

\subsection{Proof of Theorem \ref{T: Greens operator complexification}}\label{S: Greens kernel complexification}
For $y \in M$   let $G_y(x,h):=G(x,y,h)$ be the real Green's function in \eqref{analytichypoelliptic}.
 Fix  $\ep>0$, and consider the ball
 \[ B_y(\epsilon)= \{ x \in M:\; d_g(x,y) \leq \epsilon \}.\]
 Then, since $-h^2 \Delta_g +1$ is $h$-elliptic, in view of (\ref{analytic equation}),  by the semiclassical Cauchy estimates \cite[Theorem 2.6]{Jin}, for each $ x_0 \in  M \backslash B_y(\epsilon)$ there is a coordinate neighborhood  $\mathcal U \subset (M \backslash B_y(\epsilon))$  with $x_0 \in \mathcal U$  and a positive constant $C_0$, such that for all $x \in \mathcal U$ and $\alpha \in \mathbb N^n$,  
\begin{equation} \label{cauchy1}
 | \partial^{\alpha}_x G_y(x,h) | \leq C_0^{|\alpha|}  ( h^{-1} + |\alpha|)^{|\alpha|} \| G_y(\,\cdot\,, h) \|_{L^{\infty}(\mathcal U)}, 
 \end{equation}
 for all $ h \in (0,1)$.
Moreover, the estimate (\ref{cauchy1}) is locally uniform in $y \in M.$
Let $z \in {\mathcal U}^{\C}_\tau$ for a tube radius $\tau>0$ to be determined later.  By Taylor expansion around $\Re z \in \mathcal U,$ and using the Cauchy estimates (\ref{cauchy1}), we have
\begin{align} \label{cauchy2}
|G_y^{\C}(z,h)| &\leq \sum_{|\alpha| = 0}^{\infty} \frac{| \partial^{\alpha}_x G_y(x,h) |}{\alpha !} | \Im z|^{|\alpha|} \nonumber \\
& \leq  \| G_y(\,\cdot\,,h) \|_{L^{\infty}(\mathcal U)} \, \Big( \sum_{|\alpha| = 0}^{\infty}  C_0^{|\alpha|}  \frac{ ( h^{-1} + |\alpha|)^{|\alpha|}}{\alpha !} |\Im z|^{\alpha}  \Big).  
\end{align}
Let $C_1 >0$ be a large constant to be determined. Splitting the RHS of (\ref{cauchy2}) into two terms, we get
\begin{align*}
|G_y^{\C}(z,h)| &\leq T_1(z,y,h) + T_2(z,y,h),
\end{align*}
for
\begin{align}
T_1(z,y,h)&=  \| G_y(\,\cdot\,,h) \|_{L^{\infty}(\mathcal U)}\sum_{|\alpha| = 0}^{\lfloor (C_1 h)^{-1} \rfloor} \frac{ C_0^{|\alpha|} (h^{-1} + |\alpha|)^{|\alpha|} }{\alpha !} | \Im z|^{|\alpha|},  \label{E: T1} \\
T_2(z,y,h)&= \| G_y(\,\cdot\,,h) \|_{L^{\infty}(\mathcal U)} \, \sum_{|\alpha| > \lfloor (C_1 h)^{-1} \rfloor} \frac{ C_0^{|\alpha|} (h^{-1} + |\alpha|)^{|\alpha|}  }{\alpha !} | \Im z|^{|\alpha|}.   \label{E: T2}
\end{align}
To control the  term $T_1$ in \eqref{E: T1}, we use the basic estimate
$$ \sum_{|\alpha| = 0}^{\lfloor (C_1 h)^{-1} \rfloor} \frac{C_0^{|\alpha|} (h^{-1} + |\alpha|)^{|\alpha|} }{\alpha !} | \Im z|^{|\alpha|} \leq e^{C_0 (1+C_1^{-1}) |\Im z| /h},$$ and get that
$$ |T_1(z,y,h)| \leq \| G_y(\,\cdot\,,h) \|_{L^{\infty}(\mathcal U)} \, e^{ C_0(1+C_1^{-1}) |\Im z| /h}.$$
As for the second term $T_2$ in \eqref{E: T2}, we use the Stirling-type lower bounds $\alpha ! \geq (ne)^{-|\alpha|} |\alpha|^{|\alpha|}$ for all   $\alpha \in {\mathbb N}^n$  \cite[(2.9)]{Jin} 
 to get
 $$|T_2(z,y,h)| \leq  \| G_y(\,\cdot\,,h) \|_{L^{\infty}(\mathcal U)} \, \sum_{|\alpha| > h^{-1}/C_1}  |C_2 \cdot e|^{|\alpha|}| \Im z|^{|\alpha|} \leq C_3  \| G_y(\,\cdot\,,h) \|_{L^{\infty}(\mathcal U)},$$
 for some $C_3>0$, provided $z \in M_{\tau}^{\C}$ with $|\Im z| < \frac{1}{e C_2}.$
 
By the  real off-diagonal estimates in Corollary \ref{C: exponential decay real kernel},   there is a constant $C(\ep)>0$ such that
 \begin{equation} \label{realbound}
 \| G_y(\,\cdot\,,h) \|_{L^{\infty}(U)} = O(e^{-C(\ep)/h}). \end{equation}
Since $M\backslash B_y(\ep)$ is compact, the theorem follows from (\ref{realbound}) and the local bounds for $|T_1(z,y,h)|$ and $|T_2(z,y,h)|$ above, by choosing $|\Im z| < \tau(\ep)$ with $\tau(\ep) >0$ sufficiently small. Higher derivatives $|\partial_y^{\alpha} G^{\C}(z,y,h)|$ are bounded in a similar fashion.
 \qed


\subsection{Construction of  $G(x,y,h)$: Proof of Proposition \ref{P: Greens operator}}\label{S: Greens kernel}

In this section we prove Proposition \ref{P: Greens operator} by constructing an $h$-analytic parametrix $\tilde{G}(h)$ for the Green's operator $G(h)$ following closely the treatment in  \cite[Section 1]{Sj}.

Given a compact, real-analytic Riemannian manifold $(M,g)$,  consider a complex neighborhood of $T^*M$ of the form 
$$(T^*M)_\tau^\C := \left\{ (z,\zeta):\;  z \in M_\tau^\C, \;\;  |\Im {\zeta}| \leq \frac{1}{C} \langle {\zeta} \rangle  \right\}$$ 
with $C>0$ fixed sufficiently large. Here, we use the usual convention
$$  \langle \alpha_{\xi} \rangle: = \sqrt{ 1 + |\alpha_{\xi}|^2 }.$$
Following \cite{Sj}, we write $a \in S^{m,k}(T^*M)$ provided that for all $p, q \in \Z_+$
$$ \partial_{x}^{p} \partial_{\xi}^q a = O(1) h^{-m} \langle \xi \rangle^{k-|q|}$$
uniformly for $(x,\xi) \in T^*M$. We write $a \in S^{m,k}_{cl}$ if $a \sim h^{-m}(a_0 + h a_1 + ...)$ in the standard $C^{\infty}$ sense.
 The symbol $a(x,\xi,h)$ is classical analytic  (ie. $a \in S^{m,k}_{cla}$) provided $a(x,\xi,h)$ extends holomorphically to $(T^*M)_{\tau}^\C$ and the continuation (denoted by $a^\C(x,\xi,h)$) satisfies the following estimates: 
\begin{align} \label{analytic symbol}
&(i) \,\,\left|  a^\C - h^{-m} \sum_{k=0}^{ {\langle \xi \rangle}/{C_0h} } h^{k} a_k^\C \right| = O(1)e^{-\langle \xi \rangle/C_1h} \nonumber \\
&(ii) \,\, |a_j^\C(x,\xi)| \leq C_0\, C^j\,  j!\,  \langle \xi \rangle^{k-j} \nonumber \\
&(iii) \,\, \partial_{x}^{\alpha}\, \partial_{\xi}^{\beta} \,\overline{\partial_{(x,\xi)}} \,a^\C = O_{\alpha, \beta}(1) e^{-\langle \xi \rangle/C_1h} 
\end{align}
uniformly for $(x,\xi) \in (T^*M)_{\tau}^{\C}.$ In (\ref{analytic symbol}) the constant $C_0>0$ is  sufficiently large and $C_1>0$ depend on $C_0.$  In the analytic case, we henceforth write $a \sim  h^{-m} (a_0 + ha_1 +...)$ provided \eqref{analytic symbol} (i)-(iii)  hold.

Consider the phase function 
$$\phi(\alpha_x, \alpha_\xi,y) = - \left \langle  \exp_{\alpha_x}^{-1}(y)\;,\; \alpha_{\xi}\right \rangle_{\alpha_x} + \frac{i}{2}\; d_g^2(\alpha_x,y) \langle \alpha_\xi \rangle_{\alpha_x},$$
for $(\alpha_x,\alpha_\xi) \in T^*M$ and $y \in M$ with $d_g(\alpha_x,y) \leq \text{inj}(M,g)$,
and set
$$\phi^*(\alpha,x) = \overline{ \phi(\overline{\alpha},x)}.$$

\noindent Consider the cut-off function
\begin{equation}\label{E: chi}
{\rho}(x,y) := \rho (d_g(x,y)),
\end{equation}
where $ \rho \in C^{\infty}_{0}( [0,+\infty] ;[0,1])$ is a smooth cut-off function with  $\text{supp}\rho \subset [0,\text{inj} (M,g)/4]$ and  
$$ \rho(t)=1 \quad \text{for} \quad  t\in \left[0, {\text{inj}(M,g)}/{8}\right].$$
Given an elliptic symbol $ b \in S^{\frac{3n}{4},\frac{n}{4}-2}_{cla},$ consider the corresponding operator $S_{ b}(h):C^{\infty}(T^*M) \to C^{\infty}(M) $

 \begin{equation}\label{E: S}
 S_{ b}(h)u\,(x) = \int_{T^*M} e^{- \frac{i}{h} \phi^*(\alpha,x)}  \rho(\alpha_x,x)  b (\alpha,x,h) u(\alpha)   d\alpha.
 \end{equation}
Since $P(h) = -h^2 \Delta_{g} + 1$ is $h$-elliptic, there exists  $q \in S^{\frac{3n}{4},\frac{n}{4}}_{cla}$ elliptic with
$$P(h)\circ S_{ b} = S_{q}.$$
Using that $q$ is elliptic,  one can construct an $h$-analytic FBI transform  $T_a(h):C^{\infty}(M) \to C^{\infty}(T^*M)$ of the form 
\begin{equation} \label{E: T}
T_a v(\alpha;h) = \int_{M} e^{\frac{i}{h}\phi(\alpha,y)}  \rho(\alpha_x, y) a(\alpha,y;h) v(y) dv_g(y) 
\end{equation}
with $a \in S^{\frac{3n}{4},\frac{n}{4}}_{cla}$,  satisfying
\begin{equation}\label{E: R_ab}
  S_{q}(h) \circ T_a(h)  = I  + R_{ab}(h)
\end{equation}
where 
\begin{equation}\label{E: R_{ab} decay}
 | \partial_{x}^{\alpha} \partial_{y}^{\beta} R_{ab}(x,y,h) | = O_{\alpha,\beta}(e^{-C/h}).
 \end{equation}
Consequently, 
 \begin{align} \label{parametrix}
 & P(h) \circ S_{b}(h) \circ T_a(h)  = I  + R_{ab}(h). 
 \end{align}

It follows from the parametrix construction (\ref{parametrix}) that the Greens operator $G(h) = (-h^2 \Delta_{g} + 1)^{-1} \in Op_h(S^{0,-2}_{cla})$ is given by
\begin{align} \label{greens}
G(h) &=  S_{b}(h) \circ T_a(h) \circ ( I + R_{ab}(h) )^{-1} \nonumber \\
&= \tilde G(h) +\tilde R(h).\end{align}
Here, we have set 
 \begin{equation}\label{E: tilde R}
 \tilde G(h):= S_b(h) \circ T_a(h) \qquad \text{and} \qquad  \tilde R(h) =- \tilde G(h)   R_{ab}(h) (I + R_{ab}(h))^{-1}.
 \end{equation}
To compute $\tilde G(h),$ we note that from \eqref{E: S} and \eqref{E: T}, 
\begin{align}\label{E: kernel tilde G}
\tilde G(x,y; h)=\chi(x,y)
 \int_{T^*M} e^{ \frac{i}{h} \Phi(\alpha_x,\alpha_\xi, x,y)} c(\alpha_x,\alpha_\xi, x,y;h )  \; d{\alpha_\xi} d{\alpha_x},
\end{align}
where the phase function is
\begin{align}
\Phi(\alpha_x,\alpha_\xi, x,y)
&:=\phi(\alpha_x,\alpha_\xi, y)- \phi^*(\alpha_x,\alpha_\xi,x) \label{E: phase}\\
&=  \left \langle  \exp_{\alpha_x}^{-1}(x)- \exp_{\alpha_x}^{-1}(y)\;,\;\alpha_{\xi}\right \rangle_{\alpha_x} + \frac{i}{2}\; \left(d_g^2(\alpha_x,x)+d_g^2(\alpha_x,y)\right)\langle \alpha_\xi \rangle_{\alpha_x}, \notag
\end{align}
and the amplitude is
\begin{equation}\label{E: amplitude c}
c(\alpha_x,\alpha_\xi, x,y;h ):=   a(\alpha_x,\alpha_\xi,y;h) b (\alpha_x,\alpha_\xi,x,h) \rho(\alpha_x, y) \rho(\alpha_x,x) .
\end{equation}
We note that the prefactor $\chi(x,y)$ can be added in  \eqref{E: kernel tilde G} since $\chi(x,y)=1$ whenever 
$ \rho(\alpha_x, y) \rho(\alpha_x,x) \neq 0$.

Given $x,y \in M$ with $d(x,y) < \text{inj}(M,g),$  let $\alpha_x^c=\alpha_x^c(x,y)$ be the $\alpha_x$-critical point of the phase $\Phi$. We claim that 
$$\alpha_x^c(x,y) = \exp_y\Big(\frac{\exp_y^{-1}(x)}{2} \Big).$$
Indeed, in  normal coordinates  centered at $\alpha_x^c,$
 $$\partial_{\alpha_x}\Im  \Phi (\alpha_x,\alpha_\xi, x,y) \big|_{\alpha_x=\alpha_x^c}=\frac{1}{2} \partial_{\alpha_x} \left[ \left(d_g^2(\alpha_x,x)+d_g^2(\alpha_x,y)\right)\langle \alpha_\xi \rangle_{\alpha_x}\Big.\right] \Big|_{\alpha_x=\alpha_x^c} =0.$$
 Because $d_{\alpha_x}  d^2_g(\alpha_x,x)= -2 \exp_{\alpha_x}^{-1}(x)$  and  $\partial_{\alpha_x} g^{ij}(\alpha_x^c)=0$,  it follows that
$$-2(x+y)\langle \alpha_\xi \rangle_{\alpha_x^c}=0.$$
Therefore,  in normal coordinates centered at $\alpha_x^c,$ this gives $x=-y$. Thus,  the critical point $\alpha_x^c=\alpha_x^c(x,y)$ is the midpoint of the geodesic segment joining  $x$ and $y$ as claimed. \\

Since the kernel of  $\tilde G(x,y; h)$  in \eqref{E: kernel tilde G} involves an integral with amplitude supported in the set $\{\alpha_x \in M:\; d_g(y, \alpha_x) \leq 1/2\}$,  the analysis is local and 
 from now on we work in local coordinates. Then, 
$$\tilde G(x,y,h)= \chi(x,y)
 \int_{\R^n} I(\alpha_\xi, x,y,h) \, d\alpha_\xi,$$
where  
\begin{equation}\label{E: I}
I(\alpha_\xi, x,y,h):= \int_{\R^n} e^{ \frac{i}{h} \Phi(\alpha_x,\alpha_\xi, x,y)} c(\alpha_x,\alpha_\xi, x,y;h )     \; d\alpha_x.
\end{equation}

To compute $I(\alpha_{\xi},x,y,h)$ in \eqref{E: I}, we apply the method of analytic stationary phase in the $\alpha_x$-variable.

Consider the auxiliary function 
$$\Psi(\alpha_x,\alpha_\xi, x,y)=\Phi(\alpha_x,\alpha_\xi, x,y)-\Phi(\alpha_x^c,\alpha_\xi, x,y).$$
Then, $\partial_{\alpha_x}\Psi(\alpha_x^c,\alpha_\xi, x,y)=0$, $\Psi(\alpha_x^c,\alpha_\xi, x,y)=0$ with $d_{\alpha_x}^2 \Psi(\alpha_x^c,\alpha_{\xi},x,y) \sim \langle \alpha_{\xi} \rangle$ and $\Im \Psi(\alpha_x,\alpha_\xi, x,y) \geq 0$.

Let $U(x,y)\subset \C^n$ be an open neighborhood of $\alpha_x^c(x,y)$ on which the Morse Lemma holds (cf. \cite[Lemma 2.7]{Sj2} ), and set 
\begin{equation}\label{E: def of V(x,y)}
V_\R(x,y)=U(x,y)\cap \R^n. 
\end{equation}
 Further, let $\delta \in C^\infty(\R^n, \R)$ be defined as  
\begin{equation}\label{E: delta_1}
\delta(\alpha_\xi)=\frac{\delta_1}{\langle \alpha_\xi \rangle}
\end{equation}
for some $\delta_1>0$  small. For each $\alpha_\xi \in \R^n$  let $\Gamma=\Gamma(\alpha_\xi ,x,y): V_\R(x,y) \to \C^n$ be the complex contour given by $\Gamma(\alpha_\xi ,x,y) = \cup_{\alpha_x \in V_\R(x,y)} \Gamma(\alpha_x;\alpha_\xi,x,y),$ where
$$\Gamma(\alpha_x; \alpha_\xi, x,y)=\alpha_x + i \; \delta(\alpha_\xi) \; \overline{\partial_{\alpha_x}\Psi(\alpha_x,\alpha_\xi, x,y)}.$$
We choose $\delta_1$  small enough so that $\Gamma(\alpha_x;\alpha_\xi ,x,y) \subset U(x,y)$. 
With this choice of contour,
$$\Im \Psi( \Gamma(\alpha_x;\alpha_\xi, x,y),\alpha_\xi, x,y) = \delta(\alpha_\xi) |\partial_{\alpha_x}\Psi(\alpha_x,\alpha_\xi, x,y)|^2 + O(\delta^2(\alpha_\xi)  |\partial_{\alpha_x}\Psi(\alpha_x,\alpha_\xi, x,y)|^2 )$$
for all $\alpha_x \in V$.
Since $\alpha_x^c(x,y)$ is a non-degenerate critical point, 
 $$\Im \Psi(\alpha_x,\alpha_\xi, x,y) \geq C | \alpha_x -\alpha_x^c|^2$$ for some $C>0$ and all $z \in \Gamma(\alpha_\xi ,x,y)$.
 Consider the boundary surface
$S_\Gamma(\alpha_\xi ,x,y):[0,\delta(\alpha_\xi)] \times \partial V_\R(x,y) \to \C^n$  joining $V_\R(x,y)$ and ${\Gamma}(\alpha_\xi ,x,y)$ given by $ S_\Gamma(\alpha_\xi ,x,y) = \cup_{t \in [0,\delta(\alpha_\xi)]} S_\Gamma(t, \alpha_x;\alpha_\xi ,x,y),$ where
$$S_\Gamma(t, \alpha_x;\alpha_\xi ,x,y):=\alpha_x +  i \;t \; \overline{\partial_{\alpha_x}\Psi(\alpha_x,\alpha_\xi, x,y)}.$$
Let  $\Omega_\Gamma (\alpha_\xi ,x,y) \subset \C^n$  be the domain with boundary  
$$ \partial \Omega_\Gamma (\alpha_\xi ,x,y) = \Gamma (\alpha_\xi ,x,y)  \cup V_\R(x,y) \cup S_\Gamma(\alpha_\xi ,x,y).$$
  \begin{figure}[h!]
\includegraphics[height=2.6cm]{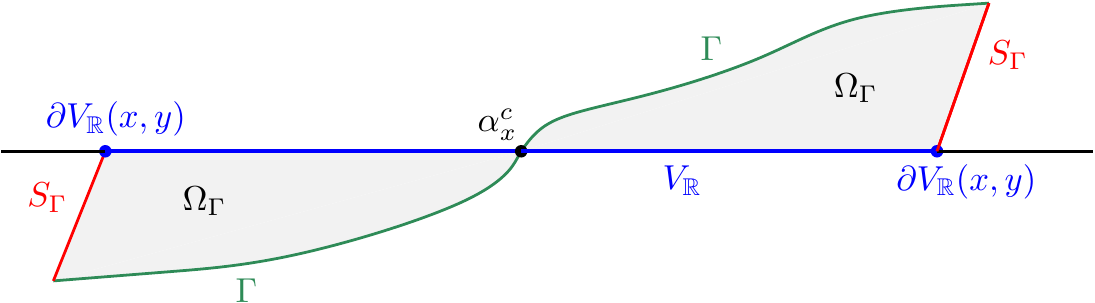}
 \end{figure}

First, from \eqref{E: I}, one can write
\begin{equation}\label{E: I 2}
I(\alpha_\xi, x,y,h)= e^{ \frac{i}{h} \Phi(\alpha_x^c,\alpha_\xi, x,y)} \int_{V_\R} e^{ \frac{i}{h} \Psi(\alpha_x,\alpha_\xi, x,y)} c(\alpha, x,y;h )     \; d\alpha_x +R_{\R^n \backslash V_{\R}}(\alpha_\xi, x,y,h),
\end{equation}
where $R_{\R^n \backslash V_{\R}}(\alpha_\xi, x,y,h)$ denotes the integral over $\R^n \backslash V_\R(x,y)$. Then, by Stoke's Theorem,
\begin{align}\label{E: I_p,ep 2} 
I(\alpha_\xi, x,y,h)&
= e^{ \frac{i}{h} \Phi(\alpha_x^c,\alpha_\xi, x,y)} \int_\Gamma e^{ \frac{i}{h} \Psi^\C(z,\alpha_\xi, x,y)} c^{\A}(z,  \alpha_\xi, x,y,h) dz  \notag\\
&\hspace{6cm} + R_\Gamma(\alpha_\xi, x,y,h) +R_{\R^n \backslash V_{\R}}(\alpha_\xi, x,y,h).
\end{align}
Here, the phase is $\Psi^\C(z,\alpha_\xi, x,y)=\Phi^\C(z,\alpha_\xi, x,y)-\Phi(\alpha_x^c,\alpha_\xi, x,y)$,  
where $\Phi^\C$ denotes the holomorphic extension of the function  $\Phi$ defined in \eqref{E: phase} in the $\alpha_x$-variable.
In the first term on the RHS of (\ref{E: I_p,ep 2}), 
$$c^{\A}(z,  \alpha_\xi, x,y,h):=a^\C(z,\alpha_\xi,y;h) b^\C (z,\alpha_\xi,x,h) \rho^{\A}(z, y) \rho^{\A}(z,x),$$
where $\rho^{\A}$ denotes the almost-analytic extension of $\rho$ (see \cite[Theorem 3.6]{Zw}).
The remainder  $ R_\Gamma$ involves integration over $\Omega_\Gamma$ and $S_\Gamma$.
An explicit description of $R_{\Gamma}, \,  R_{\R^n \backslash V_{\R}}$ and all other subsequent remainder terms  are given in \S \ref{S: remainder}.
By choosing the contour deformation space $\Omega_{\Gamma}$ sufficiently small (after possibly rescaling the parameter $t \in [0,1])$), it follows from the holomorphic Morse Lemma that there exist holomorphic local coordinates $w=( w_1, \dots, w_n)$ in a neighbourhood of $\Omega_{\Gamma}$ containing $\alpha_x^c(x,y)$ such that 
$$\Psi^\C(w, \alpha_\xi, x,y) = i \frac{(w- \alpha_x^c)^2}{2}\langle \alpha_\xi \rangle.$$

Letting $\tilde \Gamma(\alpha_\xi ,x,y)$  be the image of $\Gamma(\alpha_\xi ,x,y)$ under the change of variables $z \mapsto w,$ one can write 
\begin{equation} \label{E: I_p,ep 3}
I(\alpha_\xi, x,y,h)=  e^{ \frac{i}{h} \Phi(\alpha_x^c,\alpha_\xi, x,y)} \int_{\tilde \Gamma} e^{- \frac{(w- \alpha_x^c)^2}{2 h} \langle \alpha_\xi \rangle} {\tilde c}^{\A}(w,  \alpha_\xi, x,y,h)  dw+ (R_\Gamma+R_{\R^n \backslash V_{\R}})(\alpha_\xi, x,y,h),
\end{equation}
where $  {\tilde c}^{\A}(w,  \alpha_\xi, x,y,h):=  {c}^{\A}(w,  \alpha_\xi, x,y,h) \det\left(\frac{dz}{dw}(w, \alpha_\xi, x, y)\right).$
This choice of coordinates and the definition of $\Gamma(\alpha_\xi ,x,y)$ imply that 
$$\Re [(w- \alpha_x^c)^2 ]\geq C |w-\alpha_x^c|^2$$ for some $C>0$ and all $w \in \tilde \Gamma(\alpha_\xi ,x,y)$.

By the Implicit Function Theorem, there exists an open subset $\tilde \Gamma_0(\alpha_\xi ,x,y) \subset \tilde \Gamma(\alpha_\xi ,x,y)$ and a neighborhhod $V_{\tilde \Gamma_0} (\alpha_\xi ,x,y)\subset \R^n$ of $\alpha_x^c(x,y)$ so that 
$$\tilde \Gamma_0(\alpha_\xi ,x,y)=\{ \alpha_x + i H(\alpha_x;\alpha_\xi ,x,y): \; \alpha_x \in  V_{\tilde \Gamma_0} (\alpha_\xi ,x,y)\}$$
where $H(\,\cdot\, ;\alpha_\xi ,x,y)$ is an analytic function with $|H(\alpha_x; \alpha_\xi ,x,y)| \leq \lambda |\alpha_x|$ for some $\lambda<1$ independent of $\alpha_x$.  Then,

\begin{align} \label{E: I_p,ep 4}
I(\alpha_\xi, x,y,h)
&=  e^{ \frac{i}{h} \Phi(\alpha_x^c,\alpha_\xi, x,y)} \int_{\tilde \Gamma_0} e^{- \frac{(w- \alpha_x^c)^2}{2 h} \langle \alpha_\xi \rangle} {\tilde c}^{\A}(w,  \alpha_\xi, x,y,h)  dw \notag\\
&\hspace{3cm}+ R_\Gamma(\alpha_\xi, x,y,h)+R_{\R^n \backslash V_{\R}}(\alpha_\xi, x,y,h)+ R_{\Gamma\backslash\tilde \Gamma_0}(\alpha_\xi, x,y,h),
\end{align}
where $R_{\Gamma\backslash\tilde \Gamma_0}$ comes from replacing the domain of  the integral  in \eqref{E: I_p,ep 3} with $\Gamma(x,y)\backslash\tilde \Gamma_0(x,y)$.

As before, consider the boundary surface
$S_{\tilde \Gamma_0}:[0,1] \times \partial V_{\tilde \Gamma_0} \to \C^n$ 
$$S_{\tilde \Gamma_0}(t, \alpha_x)=\alpha_x +  i\, t \,H(\alpha_x;\alpha_\xi ,x,y)$$
joining $\tilde \Gamma_0$ with $V_{\tilde \Gamma_0}$.
Also, set $\Omega_{\tilde \Gamma_0}(\alpha_\xi ,x,y)\subset \C^n$  to be the domain whose boundary is $\tilde \Gamma_0 \cup V_{\tilde \Gamma_0}\cup S_{\tilde \Gamma_0}$. Then, another application of Stoke's Theorem   in \eqref{E: I_p,ep 4} gives
\begin{align}\label{E: I_p,ep 5}
I(\alpha_\xi, x,y,h)
&=e^{ \frac{i}{h} \Phi(\alpha_x^c,\alpha_\xi, x,y)} \int_{V_{\tilde \Gamma_0}} e^{- \frac{(\alpha_x- \alpha_x^c)^2}{2h}  \langle \alpha_\xi \rangle} {\tilde c}(\alpha_x,  \alpha_\xi, x,y,h) d\alpha_x \notag\\
&\quad+ R_\Gamma(\alpha_\xi, x,y,h)+R_{\R^n \backslash V_{\R}}(\alpha_\xi, x,y,h)+ R_{\Gamma\backslash\tilde \Gamma_0}(\alpha_\xi, x,y,h)+R_{\tilde \Gamma_0}(\alpha_\xi, x,y,h), 
\end{align}
where $R_{\tilde \Gamma_0}$ involves integration over  $\Omega_{\tilde \Gamma_0}$ and $S_{\tilde \Gamma_0}$.  

To finish the argument, we use that 
$${\tilde c}^{\A}(z,\alpha_\xi, x,y;h ):=   {\tilde a}^\C(z,\alpha_\xi,y;h) {\tilde b}^\C (z,\alpha_\xi,x,h) {\tilde \rho}^{\A}(z, y) {\tilde \rho}^{\A}(\alpha_x,x),$$ and that   $({\tilde a}^\C  \cdot  {\tilde b}^\C) (z,\alpha_\xi,x,h) $ is holomorphic in $z$. Then, by standard asymptotics for Laplace integrals  \cite[Theorem 2.1]{Sj2},
\begin{equation}\label{E: before steepest descent}
\int_{V_{\tilde \Gamma_0}}  e^{- \frac{(\alpha_x- \alpha_x^c)^2}{2h}  \langle \alpha_\xi \rangle}(a\cdot b)(\alpha_x,  \alpha_\xi, x,y,h) d\alpha_x= \frac{1}{(2\pi h)^n} \, a_G( x,y,\alpha_\xi,h) + Q_{a_G}( x,y,\alpha_\xi, h),
\end{equation}
for $a_G \in S^{0,-2}_{cl}$ defined by

\begin{equation}\label{E: a_G}
a_G( x,y,\alpha_\xi,h)=\sum_{k=0}^{\frac{1}{C_0}[\frac{1}{h}]-1} h^k w_k(x,y,\alpha_\xi),
\end{equation}
where $C_0$ is a positive constant.
Here, 
 \[w_k( x,y,\alpha_\xi)= \frac{1}{k!} \left(\frac{\Delta}{2}\right)^k (a \cdot b)(\alpha_x^c,\alpha_\xi,x,y),\]
 and 
\begin{equation}\label{E: R_{a_G}}
|Q_{a_G}( x,y,\alpha_\xi;h)|\leq C \left(1+\tfrac{1}{h}\right)^{\frac{1}{2}} e^{-\frac{1}{2h}},
\end{equation}
for some $C>0$.
Note that, in particular,    $w_0(\alpha_\xi, x,y)=(1+|\alpha_\xi|^2_{\alpha_x^c})^{-1}$.

Combining \eqref{E: I_p,ep 4} with \eqref{E: before steepest descent}  gives
\begin{align}\label{E: I final}
I(\alpha_\xi, x,y,h)
&=\frac{1}{(2\pi h)^n} e^{ \frac{i}{h} \Phi(\alpha_x^c,\alpha_\xi, x,y)} a_G( x,y,\alpha_\xi, h) + R_{\tilde G}(x,y,\alpha_\xi,h) 
\end{align}
with remainder

\begin{equation}\label{E: Q_G}
R_{\tilde G} =  R_\Gamma+R_{\R^n \backslash V_{\R}}+ R_{\Gamma\backslash\tilde \Gamma_0}+R_{\tilde \Gamma_0}+R_{1-\rho}+ R_{a_G}.
\end{equation}
In (\ref{E: Q_G}), the additional remainder term $R_{1-\rho}(\alpha_\xi,x,y,h)$ is given by
\begin{align}\label{E: R chi}
e^{ \frac{i}{h} \Phi(\alpha_x^c,\alpha_\xi, x,y)} \int_{V_{\tilde \Gamma_0}}  e^{- \frac{(\alpha_x- \alpha_x^c)^2}{2h}  \langle \alpha_\xi \rangle}(a\cdot b)(\alpha_x,  \alpha_\xi, x,y)  [1- \rho(\alpha_x, x) \rho(\alpha_x, y)]d\alpha_x,
\end{align}
and 
\begin{equation} \label{rag}
R_{a_G}(x,y,\alpha_\xi, h) = e^{ \frac{i}{h} \Phi(\alpha_x^c,\alpha_\xi, x,y)}  Q_{a_G}( x,y,\alpha_\xi,h). 
\end{equation}

Finally, 
it follows from \eqref{E: amplitude c}, \eqref{E: kernel tilde G} and \eqref{E: I final} that $\tilde G(x,y; h)$ decomposes as
\begin{equation}\label{E: G tilde}
\tilde G(x,y; h)= A_G(x,y,h)+  R_{\tilde G}(x,y,h),
\end{equation}
where we have set 
\begin{equation}\label{E: A_G}
A_G(x,y,h)=\frac{\chi(x,y) }{(2\pi h)^n} \int_{\R^n}e^{ \frac{i}{h} \Phi(\alpha_x^c,\alpha_\xi, x,y)} a_G( x,y,\alpha_\xi, h) d{\alpha_\xi},
 \end{equation}
 and
\begin{equation}\label{E: Q_G 2}
R_{\tilde G}(x,y,h)=\chi(x,y) \int_{\R^n} R_{\tilde G}(x,y,\alpha_\xi,h)d{\alpha_\xi}
\end{equation}
for $R_{\tilde G}(x,y,\alpha_\xi,h)$ defined in \eqref{E: Q_G}.\\

We now complete the proof of Proposition \ref{P: Greens operator}.

\subsubsection{Leading term $A_G(x,y,h)$} \label{S: leading term}
 Since $\alpha_x^c(x,y) = \exp_y\Big(\frac{\exp_y^{-1}(x)}{2} \Big),$ we have
 $$\Phi(\alpha_x^c(x,y),\alpha_\xi, x,y) = -2 \langle \exp_{\alpha_x^c}^{-1} (y), \alpha_\xi \rangle_{\alpha_x^c} +\frac{i}{4} d_g^2(x,y) \langle \alpha
 _\xi \rangle_{\alpha_x^c}.$$
Given $p,q \in M$, consider the parallel transport operator (along the unique shortest geodesic from $q$ to $p$) $\mathcal T_{q\to p}:T_q^*M\to T_p^*M$. This map is an isometry that satisfies
\begin{align*}
\mathcal T_{q\to p} \exp_q^{-1}(p)
&=-\exp_p^{-1}(q)\qquad \text{and}\qquad \mathcal T_{q\to p}=\mathcal T_{p\to q}^*.\label{E:Parallel Transport}
\end{align*}
Changing variables $\alpha_\xi \mapsto \eta:= \tilde{ \mathcal T}_{\alpha_x^c(x,y) \to y} (\alpha_\xi)$, where $\tilde {\mathcal  T}_{\alpha_x^c(x,y) \to y}:\R^n \to \R^n$ denotes the map induced by the choice of coordinates,  and using that $ \exp_{y}^{-1} (\alpha_x^c)= \frac{1}{2} \exp_y^{-1}(x)$, we get from \eqref{E: A_G} that in local coordinates
\begin{align*}
 A_G(x,y,h)
&=\frac{\chi(x,y) }{(2\pi h)^n}\int_{\R^n}e^{ \frac{i}{h}   \psi(x,y, \eta) }   a_G(\eta, x,y,h) d\eta,
 \end{align*}
with
 \begin{equation}\label{E: psi}
 \psi(x,y, \eta):=\langle g^{-1}_y\exp_{y}^{-1} (x),   \eta \rangle +\frac{i}{4} d_g^2(x,y) \langle g^{-1/2}_y\eta \rangle,
 \end{equation}
  and where after some abuse of notation we have set
   $$a_G (x,y,\eta,h):= a_G(x,y, \tilde{\mathcal T}_{y \to \alpha_x^c(x,y)} (\eta),h)\det\Big(\frac{d\alpha_\xi}{d\eta}(\eta, x, y)\Big)$$ 
   for $a_G \in S^{0,-2}_{cl}$ defined in \eqref{E: a_G}. In particular, since  $|\mathcal T_{y \to \alpha_x^c}(\eta)|_{\alpha_x^c}=|\eta|_y$, we have
$a_G(x,y,\eta,0)=\frac{1}{1+|g^{-1/2}_y\eta |^2}\det\big(\frac{d\alpha_\xi}{d\eta}(\eta, x, y)\big)$.
  This proves the identity \eqref{E: G theorem} for the leading term $A_{G}(x,y,h)$ in Proposition \ref{P: Greens operator}.
  

\subsubsection{Remainder  term $R_{ G}(x,y,h)$ }\label{S: remainder}
We proceed to prove statement \eqref{E: R_G theorem} in Proposition \ref{P: Greens operator}.
 In the notation of Theorem \ref{T: Greens operator complexification},
  $$ G(x,y,h)= A_G(x,y,h)+R_G(x,y,h),$$
  with
  \begin{equation*}
  R_G(x,y,h)=R_{\tilde G}(x,y,h)+ \tilde R(x,y,h).
  \end{equation*}
Here, we recall that $\tilde R(h) =- \tilde G(h)   R_{ab}(h) (I + R_{ab}(h))^{-1}$ as defined in  \eqref{E: tilde R} and
  $$R_{\tilde G}(x,y,h)= \chi(x,y) \int_{\R^n} R_{\tilde G}(\alpha_\xi,x,y,h) d{\alpha_\xi},$$
where according to \eqref{E: Q_G}
\begin{equation} \label{total remainder 1}
R_{\tilde G} =  R_\Gamma+R_{\R^n \backslash V_{\R}}+ R_{\Gamma\backslash\tilde \Gamma_0}+R_{\tilde \Gamma_0}+R_{a_G}+ R_{1-\rho}. \end{equation}
  We now prove the exponential decay in $h$ for each  of the remainder terms comprising $R_G(x,y,h)$. The exponential decay of $\partial_x^\alpha \partial_y^\beta R_G(x,y,h)$  is proved in the same way.\\

\noindent \emph{Remainders $R_{\R^n \backslash V_{\R}}$ and $R_{\Gamma\backslash\tilde \Gamma_0}$.} The term $R_{\R^n \backslash V_{\R}}$ (resp. $R_{\Gamma\backslash\tilde \Gamma_0}$) is a result of shrinking the domain of integration $\R^n$ to $V_\R(x,y) \subset \R^n$  (resp. $\Gamma$ to $\tilde \Gamma_0 \subset \Gamma$). Namely,
\begin{equation}\label{E: R_V}
R_{\R^n \backslash V_{\R}}(x,y,h)=\chi(x,y)\int_{\R^n} \int_{\R^n \backslash V_\R} e^{ \frac{i}{h} \Phi(\alpha_x,\alpha_\xi, x,y)} c(\alpha_x, \alpha_\xi, x,y;h )     \; d\alpha_x d\alpha_\xi.
\end{equation}
To study the decay of $R_{\R^n \backslash V_{\R}}$, assume without loss of generality that $V_\R(x,y)$ is a cube centered at $\alpha_x^c(x,y)$ with side length $2\delta_0$ with $\delta_0>0$ independent of $x$ and $y$:
 \begin{equation}\label{E: V_R}
V_\R(x,y)=\{\alpha_x\in \R^n:\;\; |\alpha_x^{(k)}- (\alpha_x^c(x,y))^{(k)}|<\delta_0,  \quad  k =1, \dots, n\}.
\end{equation}
 Given $\alpha_x \in  \R^n \backslash V_\R(x, y ),$  we have that  either $d_g(y, \alpha_x) > \delta_0/2$ or $d_g(x, \alpha_x)>\delta_0/2$.  Consequently, since $\Im  \Phi(\alpha_x,\alpha_\xi, x, y) =  \frac{1}{2}(d_g^2(\alpha_x,x)+d_g^2(\alpha_x,y))\langle \alpha_\xi \rangle_{\alpha_x}$, 
\[ \Im  \Phi(\alpha_x,\alpha_\xi, x, y) \geq  \frac{\delta_0^2}{8} \langle \alpha_\xi \rangle_{\alpha_x}.\]
Thus, there exists $C(\delta_0)>0$ with
 $$R_{\R^n \backslash V_{\R}} = O(e^{-C(\delta_0)/h}).$$
 The decay for $R_{\Gamma\backslash\tilde \Gamma_0}$ is proved in the same way.\\

\noindent \emph{Remainders  $ R_\Gamma$ and  $ R_{\tilde \Gamma_0}$.} 
The  term $ R_\Gamma$  (resp. $ R_{\tilde \Gamma_0}$) is the result of an application of Stoke's Theorem and consists of an integral over  $S_\Gamma$ and an integral over $\Omega_{\Gamma}$ (resp. $S_{\tilde \Gamma_0}$ and $\Omega_{\tilde \Gamma_0}$). More precisely,
$$R_\Gamma= R_{\Omega_\Gamma}+R_{S_\Gamma},$$
where
\begin{equation}\label{E: R_S}
    R_{S_\Gamma}(x, y;h)   
= \chi(x,y)\int_{\R^n} \int_{S_\Gamma} e^{ \frac{i}{h} \Phi^\C(w,\alpha_\xi, x,y)} c^{\A}(w,  \alpha_\xi, x,y) dwd\alpha_\xi,
\end{equation}
and
\begin{equation}\label{E: R_U}
 R_{\Omega_\Gamma}(x, y;h)   
=\chi(x,y)\int_{\R^n} \int_{\Omega_\Gamma} e^{ \frac{i}{h} \Phi^\C(w,\alpha_\xi, x,y)} \overline{\partial_w} c^{\A}(w,  \alpha_\xi, x,y) dwd\alpha_\xi.
\end{equation}

We first prove decay for $ R_{S_\Gamma}$. As before, for $w \in \partial V_\R(x,y),$  either $d_g(w,x)>\delta_0/2$ or $d_g(w, y)\geq \delta_0/2$. Also,  by choosing $\delta_1$ in \eqref{E: delta_1} sufficiently small in terms of $\delta_0$, one can arrange  that 
\[d_g \big(S_\Gamma(\alpha_\xi, x, y) \,,\,  \partial V_\R(x,y) \big) < \delta_0/4.\]
By Taylor expanding  $ \Phi^\C(w,\alpha_\xi, x,y)$ at $|\Im w|=0$ we get that there exists $C>0$ so that 
$\Im \Phi^\C(w,\alpha_\xi, x,y) \geq C\langle \alpha_\xi \rangle_{\Re w}$ for all $w \in S_\Gamma(\alpha_\xi, x, y)$.  This gives the exponential decay of $ R_{S_\Gamma}$.

As for the remainder $R_{\Omega_\Gamma}(x,y,h),$  one uses the fact that the amplitude in the integral for  $R_{\Omega_\Gamma}(x,y,h)$ contains the term $\overline{\partial_w } c^{\A}(w,\alpha_\xi,x,y)$
and that 
\[\overline{\partial_w } c^{\A}(w,\alpha_\xi,x,y)= (a\cdot b)^\C(w,\alpha_\xi,x,y) \cdot  \overline{\partial_w } [\rho^{\A}(w,x)\rho^{\A}(w,y)]. \]
We know that $\rho^{\A}(\Re w, x)\rho^{\A}(\Re w,y)=1$ if both $d_g(\Re w, x) < \text{inj}(M,g)/8$ and $d_g(\Re w, y) < \text{inj}(M,g)/8$ hold.  It follows that  the integrand for $R_{\Omega_\Gamma}(x,y,h)$ has its support contained in
  $$\{ w \in \Omega_{\Gamma}:\; d_g(\Re w, x) > \text{inj}(M,g)/8\,\;\; \text{or} \, \;\;d_g(\Re w, y) > \text{inj}(M,g)/8 \}.$$ The rest of the argument is the same as that for $R_{S_\Gamma}(x,y,h)$.
The analysis of the decay of $R_{\tilde \Gamma_0}$ is analogue to that of $R_\Gamma$ so we omit it.\\

\noindent \emph{Remainder $R_{1-\rho}$}. The term $R_{1-\rho}$ arises after removing the cut-off functions from the symbol $c^{\A}$ so that the result is an analytic symbol and then one can apply analytic stationary phase for quadratic phase functions. It follows from \eqref{E: R chi} that
\begin{align}\label{E: R chi 2}
&R_{1-\rho}(x,y,h) =\notag\\
&= \int_{\R^n}  \int_{V_{\tilde \Gamma_0}}   e^{ \frac{i}{h} \Phi(\alpha_x^c,\alpha_\xi, x,y)- \frac{(\alpha_x- \alpha_x^c)^2}{2h}\langle \alpha_\xi \rangle } (a\cdot b)(\alpha_x,  \alpha_\xi, x,y)  [1- \rho(\alpha_x, x) \rho(\alpha_x, y)]d\alpha_x  d\alpha_\xi.
\end{align}
The integrand of $R_{1-\rho}(x,y,h)$ is supported in the set of $(\alpha_x,x,y) \in V_{\tilde \Gamma_0} \times M \times M$ for which $d_g(\alpha_x, x)> \text{inj}(M,g)/8$ or $d_g(\alpha_x, y)> \text{inj}(M,g)/8$. Since the variable $\alpha_x$ ranges over $V_{\tilde \Gamma_0}$, we deduce that $R_{1-\rho}(x,y,h)=0$ unless $d_g(x,y)\geq C_0$ for some $C_0>0$. The rest of the argument is the same as for $R_{S_\Gamma}(x,y,h)$.\\

\noindent \emph{Remainder $R_{a_G}$}. 
We recall from (\ref{rag}) that

\begin{align}
&R_{a_G}(x,y,h) = \notag\\
&=\int_{\R^n}  \int_{ V_{\tilde \Gamma_0}} e^{ \frac{i}{h} \left(\Phi(\alpha_x^c,\alpha_\xi, x,y)+ i \frac{(\alpha_x- \alpha_x^c)^2}{2}  \langle \alpha_\xi \rangle \right)}\left[ (a\cdot b)(\alpha_x,  \alpha_\xi, x,y) - c_h(\alpha_x,  \alpha_\xi, x,y) \right] d\alpha_x  d\alpha_\xi \label{E: R_a_G  1} \\
&\qquad -\int_{\R^n} \int_{ \R^n \backslash V_{\tilde \Gamma_0}} e^{ \frac{i}{h} \left(\Phi(\alpha_x^c,\alpha_\xi, x,y)+ i \frac{(\alpha_x- \alpha_x^c)^2}{2}  \langle \alpha_\xi \rangle \right)} \, c_h(\alpha_x,  \alpha_\xi, x,y) d\alpha_x  d\alpha_\xi, \label{E: R_a_G  2}
\end{align}
where 
\begin{equation}\label{E: c_h}
c_h(\alpha_x,  \alpha_\xi, x,y):=\sum_{|\gamma |\leq \frac{1}{C_0} [\frac{1}{h}]-1}  \frac{\partial_{\alpha_x}^{\gamma}[a\cdot b](\alpha_x^c,  \alpha_\xi, x,y)}{\gamma!}  (\alpha_x - \alpha_x^c)^\gamma. 
\end{equation}

The exponential decay of \eqref{E: R_a_G  1} follows from the fact that $ a \in S^{\frac{3n}{4},\frac{n}{4}}_{cla}$ and  $ b \in S^{\frac{3n}{4},\frac{n}{4}-2}_{cla}$. Indeed, the error term $a\cdot b - c_h,$ in the Laplace integral asymptotic (see (\ref{E: a_G}) - (\ref{E: R_{a_G}})) satisfies the estimate
\[\big|(a\cdot b - c_h)(\alpha_x,  \alpha_\xi, x,y) \big| \leq  e^{-\frac{C_1}{h} \langle \alpha_\xi \rangle} \]
 with $C_1>0$.
The exponential decay of \eqref{E: R_a_G  2} is obtained in the same way  as for $R_{\R^n \backslash V_{\R}}$.\\

\noindent \emph{Remainder $\tilde R(x,y,h)$.}
Finally, we estimate  $ \tilde R(h) =- \tilde G(h)   R_{ab}(h) (I + R_{ab}(h))^{-1}$.  From \eqref{E: G tilde}, 
  \begin{align}\label{E: tilde R 1}
 \tilde R(x,y,h)
 &=\int_{M} {A_G}(x,u,h) R_{ab}(1 + R_{ab})^{-1}(u,y,h) \, du  \notag\\
&+\int_{M} {R_{\tilde G}}(x,u,h) R_{ab}(1 + R_{ab})^{-1}(u,y,h) \, du.
 \end{align}
  To deal with  the second integral in \eqref{E: tilde R 1} one simply  uses the pointwise bound
  $|R_{\tilde G}(x,u,h)| = O(e^{-C/h})$ to get that
  $$ \int_{M} {R_{\tilde G}}(x,u,h) R_{ab}(1 + R_{ab})^{-1}(u,y,h) \, du = O(e^{-C/h}).$$
To estimate the first integral in \eqref{E: tilde R 1}, we note that
$|A_G(x,u,h)|=O(1)$,  and use that the exponential decay of  $R_{ab}(x,y,h)$ in \eqref{E: R_{ab} decay} to give $\tilde R(x,y,h) = O(e^{-C/h})$ uniformly for $x,y \in M.$ \qed





\begin{thebibliography}{1}

\bibitem[BR]{BR}
J. Bourgain and Z. Rudnick.
\newblock On the nodal sets of toral eigenfunctions.
\newblock {\em Inventiones mathematicae,} 185.1 (2011): 199-237.

\bibitem[CT]{CT}
Y.~Canzani and J.~Toth.
On the local geometry of nodal sets of Laplace eigenfunctions on compact manifolds. Preprint.

\bibitem[Hel]{Hel}
B. Helffer.
\newblock Semi-classical analysis for the Schr\"odinger operator and applications.
\newblock {\em Lecture notes in mathematics},  1336 (1988).


\bibitem[HT]{HT} L. El-Hajj and J. Toth.
\newblock Intersection bounds for nodal sets of planar Neumann eigenfunctions with interior analytic curves.
\newblock {\em Journal of Differential Geometry} (2014)(to appear) arXiv:1211.3395v2.

\bibitem[HZZ]{HZZ} B. Hanin, S. Zelditch and P. Zhou.
\newblock Nodal sets of random eigenfunctions of the isotropic Harmonic Oscillator.
\newblock{\em International Mathematics Research Notices} (2014).

\bibitem[Jin]{Jin} L. Jin.
\newblock Semiclassical Cauchy estimates and applications. Preprint arXiv: 13025363.


\bibitem[LGS]{LGS}
E.~Leichtnam, F.~Golse, and M.~Stenzel.
\newblock Intrinsic microlocal analysis and inversion formulae for the heat
  equation on compact real-analytic riemannian manifolds.
\newblock {\em Annales scientifiques de l' \'Ecole normale sup\'erieure}, 29(6):669--736, 1996.

\bibitem[Mar]{Mar} A. Martinez.
\newblock An introduction to semiclassical and microlocal analysis.
\newblock {\em Springer} (2002).


\bibitem[Sj]{Sj}
J. ~Sj\"{o}strand.
\newblock Density of resonances for strictly convex analytic obstacles.
\newblock {\em Canadian Journal of  Mathematics}, 48(2):397--447, 1996.

\bibitem[Sj2]{Sj2}
J.~Sj\"{o}�strand.
\newblock Singularit\'es analytiques microlocales. 
\newblock{\em Socie�t\'{e}� math�\'{e}matique de France} Vol. 82. No. 3, 1982.

\bibitem[TZ]{TZ}
J.~Toth and S.~Zelditch.
\newblock Counting nodal lines which touch the boundary of an analytic domain.
\newblock {\em Journal of Differential Geometry}, (81):649--686.

 \bibitem[Y1]{Y1} S.T. Yau.
 \newblock Survey on partial differential equations in differential geometry, in Seminar on Differential Geometry, 3-71, Ann. of Math. Stud., 102, Princeton Univ. Press, Princeton, NJ, 1982.
 
 \bibitem[Y2]{Y2} S.T. Yau.
 \newblock Open problems in geometry, in Differential geometry: partial differential equations on manifolds (Los Angeles, CA, 1990), 1-28, Proc. Sympos. Pure Math., 54, Part 1, Amer. Math. Soc., Providence, RI, 1993.

\bibitem[Z]{Z}
S.~Zelditch.
\newblock Pluri-potential theory on Grauert tubes of real analytic Riemannian manifolds, I. Spectral geometry. {\em   Proceedings of Symposia in Pure Mathematics,} (84):299-339.  American Mathematical  Society, Providence, RI.

\bibitem[Z2]{Z2}
S.~Zelditch.
\newblock Eigenfunctions and nodal sets. (2012). Preprint arXiv:1205.2812.


\bibitem[Zw]{Zw} M.~Zworski.
\newblock Semiclassical analysis.
\newblock {\em  American Mathematical Society}, Vol. 138. 2012.

\end{thebibliography}
\end{document}